\documentclass[12pt,a4paper,twoside]{article}
\usepackage[a4paper,margin=2.5cm,bottom=4cm]{geometry}			% Set paper to letter and set margins

\RequirePackage{fix-cm}

\usepackage{graphicx}
\usepackage{amsfonts}
\usepackage{hyperref}
\usepackage{epstopdf}
\usepackage{algorithmic}
\usepackage{tikz}
\usepackage{amsmath} % for math
\usepackage{amssymb} % for more math
\usepackage{mathabx} % for even more math
\usepackage{mathtools}
\usepackage[caption=false]{subfig}
\usepackage{multirow}
\usepackage[version=4]{mhchem}
\usepackage{overpic}

\DeclareMathOperator{\secant}{sct}
\DeclareMathOperator{\argmax}{argmax}
\DeclareMathOperator\erf{erf}
\newcommand{\E}{\mbox{I\negthinspace E}}

\usepackage[utf8]{inputenc}
\usepackage[english]{babel}

\begin{document}

\title{Global optimization of Gaussian processes \thanks{This work was supported by the Deutsche Forschungsgemeinschaft (DFG, German Research Foundation) under Germany's Excellence Strategy - Cluster of Excellence 2186 ``The Fuel Science Center'' and the project ``Improved McCormick Relaxations for the efficient Global Optimization in the Space of Degrees of Freedom'' MI 1851/4-1.}}

\author{Artur M. Schweidtmann$^1$ 
	\and Dominik Bongartz$^1$
	\and Daniel Grothe$^1$
	\and Tim Kerkenhoff$^1$
	\and Xiaopeng Lin$^1$
	\and Jaromi{\l} Najman$^1$
	\and Alexander Mitsos$^{1,2,3}$
}

\date{%
	$^1$ Process Systems Engineering (AVT.SVT), RWTH Aachen University, Aachen, Germany.\\%
	$^2$              JARA-CSD, 52056 Aachen, Germany.\\
	$^3$ Institute of Energy and Climate Research, Energy Systems Engineering (IEK-10), Forschungszentrum Jülich GmbH, 52425 Jülich, Germany. \\
	({amitsos@alum.mit.edu}).\\[2ex]%
	\today
}

\maketitle

\begin{abstract}
Gaussian processes~(Kriging) are interpolating data-driven models that are frequently applied in various disciplines.
Often, Gaussian processes are trained on datasets and are subsequently embedded as surrogate models in optimization problems.
These optimization problems are nonconvex and global optimization is desired.
However, previous literature observed computational burdens limiting deterministic global optimization to Gaussian processes trained on few data points.
We propose a reduced-space formulation for deterministic global optimization with trained Gaussian processes embedded.
For optimization, the branch-and-bound solver branches only on the degrees of freedom and McCormick relaxations are propagated through explicit Gaussian process models.
The approach also leads to significantly smaller and computationally cheaper subproblems for lower and upper bounding.
To further accelerate convergence, we derive envelopes of common covariance functions for GPs and tight relaxations of acquisition functions used in Bayesian optimization including expected improvement, probability of improvement, and lower confidence bound.
In total, we reduce computational time by orders of magnitude compared to state-of-the-art methods, thus overcoming previous computational burdens.
We demonstrate the performance and scaling of the proposed method and apply it to Bayesian optimization with global optimization of the acquisition function and chance-constrained programming.
The Gaussian process models, acquisition functions, and training scripts are available open-source within the ``MeLOn - \textbf{M}achin\textbf{e} \textbf{L}earning Models for \textbf{O}ptimizatio\textbf{n}'' toolbox~(\url{https://git.rwth-aachen.de/avt.svt/public/MeLOn}). 
\end{abstract}

\section{Introduction}
\label{sec:Introduction}
A Gaussian process~(GP) is a stochastic process where any finite collection of random variables has a multivariate Gaussian distribution; 
they can be understood as an infinite-dimensional generalization of multivariate Gaussian distributions~\cite{rasmussen2004gaussian}.
This means that predictions of GPs are Gaussian distributions that provide not only an estimate but also a variance.
GPs originate from geostatistics~\cite{krige1951statistical} and gained popularity for the design and analysis of computer experiments~(DACE) since 1989~\cite{sacks1989design}.
Furthermore, GPs are commonly applied as interpolating surrogate models across various disciplines including 
biotechnology~\cite{freier2016framework,ulmasov2016bayesian,mehrian2018maximizing,bradford2018b,del2019review},
chemical engineering~\cite{mcbride2019overview,Caballero.2008,davis2007kriging,davis2008kriging,davis2010centroid,eason2016trust,Helmdach.2017}, 
chemistry~\cite{schweidtmann2018machine,amar2019machine}, 
and deep-learning~\cite{snoek2012practical}.
Note that GP regression is also often referred to as Kriging.
In many applications, GPs are trained on a data set and are subsequently embedded in an optimization problem, e.g., to identify an optimal operating point of a process.
Moreover, many derivative-free solvers for expensive-to-evaluate black-box functions actually train GPs and optimize their predictions (e.g., Bayesian optimization algorithms~\cite{Shahriari.2016,jones1998efficient,ulmasov2016bayesian,bradford2018efficient} and other adaptive sampling approaches~\cite{davis2007kriging,davis2008kriging,davis2010centroid,Cozad.2014,eason2016trust,boukouvala2016global,boukouvala2017argonaut}).
In Bayesian optimization, the optimum of an acquisition function determines the next sampling point~\cite{Shahriari.2016}.
The vast majority of these optimizations have been performed by local solution approaches~\cite{kim2019local,wilson2018maximizing} and a few by stochastic global~\cite{bradford2018efficient}.
Our contribution focuses on the deterministic global solution of optimization problems with trained GPs embedded and on applications in process systems engineering. \newline \indent
GPs are commonly used to learn the input-output behavior of unit operations~\cite{caballero2008rigorous,Caballero.2008,quirante2015optimization,quirante2015rigorous,lin2017process,kessler2019globalCES,kessler2019globalAS}, complete flowsheets~\cite{hasan2012modeling}, or thermodynamic property relations from data~\cite{mcbride2017integrated}.
Subsequently, the trained GPs are often combined with nonlinear mechanistic process models leading to hybrid mechanistic / data-driven models~\cite{von2014hybrid,kahrs2007validity,glassey2018hybrid,mogk2002application} which are optimized.
Most of the previous works on optimization with GPs embedded rely on local optimization techniques.
Caballero and Grossmann, for instance, train GPs on data obtained from a rigorous divided wall column simulation.
Then, they iteratively optimize the operation of the column~(modeled by GPs) using SNOPT~\cite{gill2005snopt}, sample new data at the solution point, and update the GPs~\cite{caballero2008rigorous,Caballero.2008}.
Later, Caballero and co-workers extend this work to distillation sequence superstructure problems~\cite{quirante2015optimization,quirante2015rigorous}.
In~\cite{quirante2015rigorous}, the authors solve the resulting mixed-integer nonlinear programs~(MINLPs) using a local solver in GAMS.
Therein, the GP estimate is computed via an external function in Matlab which leads to a reduced optimization problem size visible to the local solver in GAMS.
However, all these local methods have the drawback that they can lead to suboptimal solutions, because the resulting optimization problems are nonconvex. 
This nonconvexity is induced by the covariance functions of the GPs as well as often the mechanistic part of the hybrid models.  \newline \indent
Deterministic global optimization can guarantee to identify globally optimal solutions within finite time to a given nonzero tolerance~\cite{Horst.1996}. 
In a few previous studies, deterministic global optimization with GPs embedded was done using general-purpose global solvers. 
For instance, in the black-box optimization algorithms ALAMO~\cite{Cozad.2014} and ARGONAUT~\cite{boukouvala2017argonaut}, GPs are included as surrogate models and are optimized using BARON~\cite{Tawarmalani.2005} and ANTIGONE~\cite{Misener.2014}, respectively.
However, computational burdens were observed that limit applicability, e.g., in terms of the number of training points.
Cozad et al.~\cite{Cozad.2014} state that GPs are accurate but ``difficult to solve using provable derivative-based optimization software''.
Similarly, Boukouvala and Floudas~\cite{boukouvala2017argonaut} state that the computational cost becomes a limiting factor because the number of nonlinear terms of GPs equals the product of the number of interpolated points~($N$) and the dimensionality~($D$) of the input domain.
More recently, Ke{\ss}ler et al.~\cite{kessler2019globalCES,kessler2019globalAS} optimized the design of nonideal distillation columns by a trust-region approach with GPs embedded.
Therein, optimization problems with GPs embedded are solved globally using BARON~\cite{Tawarmalani.2005} within relatively long CPU times ($10^2 - 10^5$ CPU seconds on a personal computer). \newline \indent
As mentioned earlier, Quirante et al.~\cite{quirante2015rigorous} call an external Matlab function to compute GP estimates with a local solver in GAMS.
As an alternative approach, they also solve the problem globally using BARON in GAMS by providing the full set of GP equations as equality constraints.
This leads to additional intermediate optimization variables besides the degrees of freedom of the problem. 
Similar to other studies, they observe that their formulation is only practical for a small number of GP surrogates and training data points to avoid large numbers of variables and constraints~\cite{quirante2015rigorous}. 
We refer to the problem formulation where the GP is described by equality constraints and additional optimization variables as a full-space~(FS) formulation.
It is commonly used in modeling environments, e.g., GAMS, that interface with state-of-the-art global solvers such as ANTIGONE~\cite{Misener.2014}, BARON~\cite{Tawarmalani.2005}, and  SCIP~\cite{Maher.2017}.  \newline \indent
An alternative to the FS is a reduced-space~(RS) formulation where some optimization variables are eliminated using explicit constraints.
This reduced problem size leads to a lower number of variables for branching as well as potentially smaller subproblems.
The former has some similarity to selective branching~\cite{Epperly.1997} (c.f. discussion in~\cite{bongartz2020_diss}).
The exact size of the subproblems for lower bounding and bound tightening depends on the method for constructing relaxations.
In particular, when constructing relaxations in the RS using McCormick~\cite{McCormick.1976}, alphaBB~\cite{androulakis1995alphabb} or natural interval extensions, the resulting lower bounding problems are much smaller compared to the auxiliary variable method~(AVM)~\cite{smith1997,Tawarmalani.2002}.
Therefore, any global solver can in principle handle RS but some methods for constructing relaxations appear more promising to benefit from the RS~\cite{bongartz2020_diss}.
We have recently released the open-source global solver MAiNGO~\cite{MAiNGO} which uses the MC\texttt{++} library~\cite{Chachuat.2015} for automatic propagation of McCormick relaxations through computer code~\cite{Mitsos.2009}.
We have shown that the RS formulation can be advantageous for 
flowsheet optimization problems~\cite{Bongartz.2017,bongartz2019deterministic} and problems with artificial neural networks embedded~\cite{Schweidtmann2019detGlobalANN,rall2019rational,hullen2019managing}. 
In the context of Bayesian optimization, Jones et at.~\cite{jones1998efficient} develop valid overestimators of the expected improvement (EI) acquisition function in the RS.
However, their relaxations rely on interval extensions and optimization-based relaxations limited to a specific covariance function; they do not derive envelopes.
Furthermore, they do not provide convex underestimators which are in general necessary to embed GPs in optimization problems. \newline \indent
In this work, we show that the proposed RS outperforms a FS formulation for problems with GPs embedded by speedup factors of several magnitudes.
Moreover, this speedup increases with the number of training points.
This is mainly due to the fact that the number of variables and constraints in the FS scales linearly with the number of data points $\mathcal{O}(N+D)$ while in the RS it depends only on the input dimensionality, i.e., $\mathcal{O}(D)$.
To further accelerate convergence, we derive envelopes of covariance functions for GPs and tight relaxations of acquisition functions, which are commonly used in Bayesian optimization. 
Finally, we solve a chance-constrained optimization problem with GPs embedded and we perform global optimization of an acquisition function.
The GP training methods and models are provided as an open-source toolbox called ``MeLOn - \textbf{M}achin\textbf{e} \textbf{L}earning Models for \textbf{O}ptimizatio\textbf{n}'' under the Eclipse public license~\cite{MeLOn_Git}.
The resulting optimization problems are solved using our global solver MAiNGO~\cite{MAiNGO}.
Note that the MeLOn toolbox is also automatically included as a submodule in our new MAiNGO release.
\section{Gaussian Processes}
\label{sec:GP_Opt_Background}
In this section, GPs are briefly introduced (c.f.~\cite{rasmussen2004gaussian}). 
We first describe the GP prior distribution, i.e., the probability distribution before any data is taken into account. 
Then, we describe the posterior distribution, which results from conditioning the prior on training data. 
Finally, we describe how hyperparameters of the GP can be adapted to data by a maximum a posteriori (MAP) estimate. 
\subsection{Prior}
\label{sec:GP_Opt_Background_prior}
A GP prior is fully described by its mean function $m(\boldsymbol{x})$ and positive semi-definite covariance function $k(\boldsymbol{x},\boldsymbol{x}')$~(also known as kernel function).
We consider a noisy observation $y$ from a function $\tilde{f}(\boldsymbol{x})$ with  $y(\boldsymbol{x}) \coloneqq \tilde{f}(\boldsymbol{x}) + \varepsilon_{\mathrm{noise}}$, whereby the output noise $\varepsilon_{\mathrm{noise}}$ is  independent and identically distributed (i.i.d.) with $\varepsilon_{\mathrm{noise}} \sim \mathcal{N}(0,{\sigma}_{\mathrm{noise}}^2 )$. 
We say $y$ is distributed as a GP, i.e., $y \sim \mathcal{GP}(m(\boldsymbol{x}),k(\boldsymbol{x}, \boldsymbol{x}'))$ with 
\begin{align}
m(\boldsymbol{x}) &\coloneqq \E\big[\tilde{f}(\boldsymbol{x})\big], \notag \\ 
k(\boldsymbol{x}, \boldsymbol{x}') &\coloneqq \E \big[~( \medspace y(\boldsymbol{x}) - m(\boldsymbol{x}) \medspace )~( \medspace y(\boldsymbol{x}') - m(\boldsymbol{x}) \medspace )^{\mathrm{T}} \big]. \notag
\end{align}
Without loss of generality, we assume that the prior mean function is $m(\boldsymbol{x})=0$.
This implies that we train the GP on scaled data such that the mean of the training outputs is zero.
A common class of covariance functions is the Mat\'ern class. 
\begin{equation}
k_{\text{Mat\'ern}}(\boldsymbol{x}, \boldsymbol{x}') \coloneqq \sigma_f^2 \frac{2^{1-\nu}}{\Gamma(\nu)} \left( \sqrt{2 \nu} r \right)^\nu K_\nu \left( \sqrt{2 \nu} r \right), \notag
\end{equation}
where 
$\sigma_f^2$ is the output variance,
$r \coloneqq \sqrt{(\boldsymbol{x} - \boldsymbol{x}')^{\mathrm{T}} \boldsymbol{\Lambda}~(\boldsymbol{x} - \boldsymbol{x}')}$ is a weighted Euclidean distance, 
$\boldsymbol{\Lambda} \coloneqq \mathrm{diag}(\lambda_1^2, \cdots , \lambda_i^2, \cdots \lambda_{n_x}^2)$ is a length-scale matrix with $\lambda_i \in \mathbb{R}$, 
$\Gamma(\cdot)$ is the gamma function, 
and $K_\nu(\cdot)$ is the modified Bessel function of the second kind.
The smoothness of Mat\'ern covariance functions can be adjusted by the positive parameter $\nu$. 
When $\nu$ is a half-integer value, the Mat\'ern covariance function becomes a product of a polynomial and an exponential~\cite{rasmussen2004gaussian}. 
Common values for $\nu$ are $1/2,~3/2,~5/2,$ and $\infty$, i.e., the most widely-used squared exponential covariance function.
We derive envelopes of these covariance functions in Section~\ref{sec:Relaxations_of_Covariance_Functions} and implement them within MeLOn~\cite{MeLOn_Git}.
Also, a noise term, $\sigma_{\mathrm{n}}^2 \cdot \delta(\boldsymbol{x},\boldsymbol{x}')$, can be added to any covariance function where $\sigma_{\mathrm{n}}^2$ is the noise variance and $\delta(\boldsymbol{x},\boldsymbol{x}')$ is the Kronecker delta function.
The hyperparameters of the covariance function are adjusted during training and are jointly noted as $\boldsymbol{\theta} =$ $ [\lambda_1,..,\lambda_d,\sigma_f, \sigma_{\mathrm{n}}]$.
Herein, a log-transformation is common to prevent negative values during training. 
\subsection{Posterior}
\label{sec:GP_Opt_Background_posterior}
The GP posterior is obtained by conditioning the prior on observations.
We consider a set of $N$ training inputs $\boldsymbol{\mathcal{X}}=\{\boldsymbol{x}_1^{(\boldsymbol{\mathcal{D}})}, ...,\boldsymbol{x}_N^{(\boldsymbol{\mathcal{D}})} \}$ where $\boldsymbol{x}_{i}^{(\boldsymbol{\mathcal{D}})} = [x_{i,1}^{(\boldsymbol{\mathcal{D}})}, ...,x_{i,D}^{(\boldsymbol{\mathcal{D}})}]^T$ is a $D$-dimensional vector.
Note that we use the superscript ${(\boldsymbol{\mathcal{D}})}$ to denote the training data.
The corresponding set of scalar observations is given by $\boldsymbol{\mathcal{Y}}=\{y_1^{(\boldsymbol{\mathcal{D}})}, ...,y_N^{(\boldsymbol{\mathcal{D}})} \}$.
Furthermore, we define the vector of scalar observations $\boldsymbol{y} =$ 
$[y_1^{(\boldsymbol{\mathcal{D}})}, ...,y_N^{(\boldsymbol{\mathcal{D}})}]^T$ $\in \mathbb{R}^{N}$.
The posterior GP is obtained by \textsc{Bayes}' theorem: 
\begin{equation}
\tilde{f}(\boldsymbol{x}) \sim \mathcal{GP}(m(\boldsymbol{x}), k(\boldsymbol{x},\boldsymbol{x}') | \boldsymbol{\mathcal{X}},\boldsymbol{\mathcal{Y}}) = \mathcal{N} \left(
m_{\boldsymbol{\mathcal{D}}}( \boldsymbol{x}),
k_{\boldsymbol{\mathcal{D}}}( \boldsymbol{x},\boldsymbol{x}')
\right) \notag
\end{equation}
with 
\begin{align}
m_{\boldsymbol{\mathcal{D}}}( \boldsymbol{x}) &= \boldsymbol{K}_{\boldsymbol{x},\mathcal{X}} \left( \boldsymbol{K}_{\mathcal{X},\mathcal{X}} \right)^{-1}  \boldsymbol{y},  \label{eq:GP_prediction}\\
k_{\boldsymbol{\mathcal{D}}}(\boldsymbol{x}) &= {K}_{\boldsymbol{x},\boldsymbol{x}} - \boldsymbol{K}_{\boldsymbol{x},\mathcal{X}} \left( \boldsymbol{K}_{\mathcal{X},\mathcal{X}} \right)^{-1} \boldsymbol{K}_{ \mathcal{X},\boldsymbol{x} }, \label{eq:GP_variance}
\end{align}
where the covariance matrix of the training data is given by $\boldsymbol{K}_{\mathcal{X},\mathcal{X}} \coloneqq \left[ k(\boldsymbol{x}_i,\boldsymbol{x}_j) \right] \in \mathbb{R}^{N \times N}$, the covariance vector between the candidate point $\boldsymbol{x}$ and the training data is given by
$\boldsymbol{K}_{\boldsymbol{x},\mathcal{X}} \coloneqq \left[ k(\boldsymbol{x},\boldsymbol{x}_1^{(\boldsymbol{\mathcal{D}})}), ..,k(\boldsymbol{x},\boldsymbol{x}_N^{(\boldsymbol{\mathcal{D}})}) \right] \in \mathbb{R}^{1 \times N}$, $\boldsymbol{K}_{ \mathcal{X},\boldsymbol{x} } = {\boldsymbol{K}_{ \boldsymbol{x},\mathcal{X} }}^T$, and ${K}_{\boldsymbol{x},\boldsymbol{x}} \coloneqq k(\boldsymbol{x},\boldsymbol{x})$.
Equations \eqref{eq:GP_prediction}, \eqref{eq:GP_variance} describe essentially the predictions of a GP and are implemented within MeLOn.
\subsection{Maximum A Posteriori} 
\label{sec:GP_Opt_Background_MAP}
In order to find appropriate hyperparameters  $\boldsymbol{\theta}$ for a given problem, we use a MAP estimate which is known to be advantageous compared to the maximum likelihood estimation~(MLE) on small data sets~\cite{sundararajan2000predictive}.
Using the MAP estimate, the hyperparameters are identified by maximizing the probability that the GP fits the training data, i.e., $\boldsymbol{\theta}_{\mathrm{opt}} \coloneqq \argmax\limits_{\boldsymbol{\theta}}\mathcal{P}\left( \boldsymbol{\theta} |\boldsymbol{\mathcal{X}}, \boldsymbol{\mathcal{Y}} \right)$.
Analytical expressions for $\mathcal{P}\left( \boldsymbol{\theta} |\boldsymbol{\mathcal{X}}, \boldsymbol{\mathcal{Y}} \right)$ and its derivatives w.r.t. the hyperparameters can be found in the literature~\cite{rasmussen2004gaussian}.
We provide a Matlab training script in MeLOn that is based on our previous work~\cite{bradford2018efficient}.
Therein, we assume an independent Gaussian distribution as a prior distribution on the log-transformed hyperparameters, i.e., $\theta_i \sim \mathcal{N} \left(\mu_i, \sigma_i^2 \right)$.
\section{Optimization Problem Formulations}
\label{sec:GP_Opt_Method}
In the simplest and common in the literature case, the (scaled) inputs of a GP are degrees of freedom of the optimization problem with $\boldsymbol{x} \in \tilde{X} = [\boldsymbol{x}^{L},\boldsymbol{x}^{U}]$.
For given $\boldsymbol{x}$, the dependent (or intermediate) variables $\boldsymbol{z}$ can be computed by the solution of $\boldsymbol{h}(\boldsymbol{x},\boldsymbol{z})=\boldsymbol{0}, \quad \boldsymbol{h}:\tilde{X} \times \mathbb{R}^{n_z} \to \mathbb{R}^{n_z}$.
This corresponds to Equations~\eqref{eq:GP_prediction} and~\eqref{eq:GP_variance} where the estimate ($m_{\boldsymbol{\mathcal{D}}}$) and variance ($k_{\boldsymbol{\mathcal{D}}}$) of the GP are computed respectively. 
Note that in this case, we can solve explicitly for $m_{\boldsymbol{\mathcal{D}}}$ and $k_{\boldsymbol{\mathcal{D}}}$.\newline \indent
The realization of the objective function $f$ depends on the application.
In many applications, it depends on the estimate of the GP, i.e.,  $f(m_{\boldsymbol{\mathcal{D}}})$~(c.f. Subsection~\ref{subsec:Numerical Example 1}).
In Bayesian optimization, the objective function is called acquisition function and usually depends on the estimate and variance of the GP, i.e., $f(m_{\boldsymbol{\mathcal{D}}},k_{\boldsymbol{\mathcal{D}}})$~(c.f. Subsection~\ref{subsec:Numerical Example 3}). 
Finally, additional constraints might depend on the inputs of the GP, its estimate, and variance, i.e., $\boldsymbol{g}(\boldsymbol{x},m_{\boldsymbol{\mathcal{D}}},k_{\boldsymbol{\mathcal{D}}}) \le \boldsymbol{0}$.
In more complex cases, multiple GPs can be combined in one optimization problem (c.f. Subsection~\ref{subsec:Numerical Example 2}). \newline \indent
In the following, we describe two optimization problem formulations for problems with trained GPs embedded: the commonly used FS formulation in Subsection~\ref{subsec:GP_Full-space formulation} and the RS formulation in Subsection~\ref{subsec:GP_Reduced-space formulation}.
Both problem formulations are equivalent reformulations in the sense that they have the same optimal solution.
However, the formulation significantly affects problem size and performance of global optimization solvers.
\subsection{Full-Space formulation}
\label{subsec:GP_Full-space formulation}
In the FS formulation, the nonlinear equations $\boldsymbol{h}(\boldsymbol{x},\boldsymbol{z})=\boldsymbol{0}$ are provided as equality constraints and the intermediate dependent variables $\boldsymbol{z} \in Z$ are optimization variables.
A general FS problem formulation is:
\begin{align}
& \underset{\boldsymbol{x} \in \tilde{X}, \boldsymbol{z} \in Z}{\min} 	&& f(\boldsymbol{x},\boldsymbol{z}) 		 \tag*{(FS)}	\\
& \text{s.t.} 													&& \boldsymbol{h}(\boldsymbol{x},\boldsymbol{z}) = \boldsymbol{0}, \qquad \boldsymbol{g}(\boldsymbol{x},\boldsymbol{z}) \le \boldsymbol{0} \notag
\end{align}
In general, there exist multiple valid FS formulations for optimization problems.
In Section~3 of the electronic supplementary information (ESI), we provide a representative FS formulation for the case where the estimate of a GP is minimized.
This is also the FS formulation that we use in our numerical examples (c.f., Section~\ref{subsec:Numerical Example 1}).
\subsection{Reduced-Space Formulation}
\label{subsec:GP_Reduced-space formulation}
In the RS formulation, the equality constraints are solved for the intermediate variables and substituted in the optimization problem~(c.f.~\cite{Bongartz.2017}).
A general RS problem formulation in the context of optimization with a GP embedded is:
\begin{align}
& \underset{\boldsymbol{x} \in \tilde{X}}{\min} 					&& f(m_{\boldsymbol{\mathcal{D}}}(\boldsymbol{x}),k_{\boldsymbol{\mathcal{D}}}(\boldsymbol{x}))		 		&	\tag*{(RS)}		\\
& \text{s.t.} 											&& \boldsymbol{g}(\boldsymbol{x},m_{\boldsymbol{\mathcal{D}}}(\boldsymbol{x}),k_{\boldsymbol{\mathcal{D}}}(\boldsymbol{x})) \le \boldsymbol{0} & 	\notag
\end{align}
Herein, the B\&B solver operates only on the degrees of freedom $\boldsymbol{x}$ and no equality constraints are visible to the solver.
In GPs, the estimate and variance are explicit functions of the input~(Equations~\eqref{eq:GP_prediction},~\eqref{eq:GP_variance}).
Thus, we can directly formulate a RS formulation.
The RS formulation effectively combines those equations and hides them from the B\&B algorithm. 
This results in a total number of $D$ optimization variables, zero equality constraints, and no additional optimization variables $\boldsymbol{z}$.
Thus, the RS formulation requires only bounds on $\boldsymbol{x}$. \newline \indent
Note that the direct substitution of all equality constraints is not always possible when multiple GPs are combined with mechanistic models, e.g., in the presence of recycle streams.
Here, a small number of additional optimization variables and corresponding equality constraints can remain in the RS formulation~\cite{Bongartz.2017}.
As an alternative, relaxations for implicit functions can also be derived~\cite{Stuber.2015,Wechsung.2015}.
Moreover, we have previously observed that a hybrid between RS and FS formulation can be optimal for some optimization problems~\cite{bongartz2019deterministic}. 
In this work, we compare the RS and the FS formulation and do not consider any hybrid problem formulations.
\section{Convex and Concave Relaxations}
\label{sec:Relaxations_of_Relevant_Functions_for_GPs_and_Bayesian_Optimization}
The construction of relaxations, i.e., convex function underestimators ($F^{cv}$) and concave function overestimators ($F^{cc}$), is essential for B\&B algorithms.
In our open-source solver MAiNGO, we use the (multivariate) McCormick method~\cite{McCormick.1976,Tsoukalas.2014} to propagate relaxations and their subgradients~\cite{Mitsos.2009} through explicit functions using the MC\texttt{++} library~\cite{Chachuat.2015}.
However, the McCormick method often does not provide the tightest possible relaxations, i.e., the envelopes.
In this section, we derive tight relaxations or envelopes of functions that are relevant for GPs and Bayesian optimization.
The functions and their relaxations are implemented in MC\texttt{++}.
When using these intrinsic functions and their relaxations in MAiNGO, the (multivariate) McCormick method is only used for the remaining parts of the model.
Note that the derived relaxations are used within MAiNGO while BARON does not allow for implementation of custom relaxations or piecewise defined functions.
\subsection{Covariance Functions}
\label{sec:Relaxations_of_Covariance_Functions}
The covariance function is a key element of GPs that occurs $N$ times within the optimization of estimate or variance.
Thus, tight relaxations are highly desirable.
In this subsection, we derive envelopes for common Mat\'ern covariance functions.
We consider univariate covariance functions, i.e., $k_{\nu}: \mathbb{R} \to \mathbb{R}$, with input $d = (\boldsymbol{x} - \boldsymbol{x}')^{\mathrm{T}} \boldsymbol{\Lambda}~(\boldsymbol{x} - \boldsymbol{x}')\ge 0$.
This is possible because we consider stationary covariance functions that are invariant to translations in the input space. 
Common Mat\'ern covariance functions use $\nu=1/2,~3/2,~5/2$ and $\infty$ and are given by:
\begin{align}
&k_{\nu=1/2}(d) \coloneqq \exp \left( -\sqrt{d} \right), \qquad
k_{\nu=3/2}(d) \coloneqq \left( 1 + \sqrt{3}~\sqrt{d}  \right) \cdot \exp \left(- \sqrt{3}~ \sqrt{d}  \right) \notag \\
&k_{\nu=5/2}(d) \coloneqq \left( 1 + \sqrt{5}~\sqrt{d}  +\frac{5}{3}~d  \right) \cdot \exp \left( -\sqrt{5}\sqrt{d}  \right), \qquad
k_{SE}(d) \coloneqq \exp \left( -\frac{1}{2}~d \right), \notag
\end{align}
where $k_{SE}$ is the squared exponential covariance function with $\nu \rightarrow \infty$.
We find that these four covariance functions are convex because their Hessian is positive semidefinite.
Thus, the convex envelope is given by $F^{cv}(d) = k(d)$ and the concave envelope by the secant $F^{cc}(d) = \secant(d)$ where 
$\secant(d) =$ $\frac{k(d^U) - k(d^L)}{d^U - d^L} d + \frac{d^U k(d^L) - d^L k(d^U) }{d^U - d^L}$ on a given interval $[d^L, d^U]$.
As the McCormick composition and product theorems provide weak relaxations of $k_{\nu=3/2}$ and $k_{\nu=5/2}$ (c.f. ESI Section~1), we implement these functions and their envelopes in our library of intrinsic functions in MC\texttt{++}.
Furthermore, natural interval extensions are not exact for $k_{\nu=3/2}$ and $k_{\nu=5/2}$.
Thus, we also provide exact interval bounds based on the monotonicity. \newline   \indent
It should be noted that covariance functions are commonly given as a function of the weighted Euclidean distance $r=\sqrt{d}$.
However, we chose to use d instead for three main reasons:
(1) $\boldsymbol{x}$ is usually a degree of freedom.
Thus, the computation of $r$ would lead to potentially weaker relaxations for $k_{\nu=5/2}$ and $k_{SE}$.
(2) The derivative of $k_{\nu=3/2}(\cdot)$, $k_{\nu=5/2}(\cdot)$, and $k_{SE}(\cdot)$ is defined at $d=0$ while the derivative of the square root function is not.
(3) The covariance functions $\hat{k}_{v=3/2}: r \mapsto k_{v=3/2}(r^2)$, $\hat{k}_{v=5/2}: r \mapsto k_{v=5/2}(r^2)$, and $\hat{k}_{SE}: r \mapsto k_{SE}(r^2)$ are nonconvex in $r$, so deriving the envelopes would be nontrivial. \newline \indent
Finally, it can be noted that we did not derive envelopes of $k_{\text{Mat\'ern}}(\boldsymbol{x}, \boldsymbol{x}')$, because the variable input dimensions poses difficulties in implementation and the multidimensionality is a challenge for derivation of envelopes.
Nevertheless, the McCormick composition theorem applied to $k_{\nu}(d(\boldsymbol{x},\boldsymbol{x}'))$ yields relaxations that are exact at the minimum of $k_{\text{Mat\'ern}}$ because the natural interval extensions of the weighted squared distance $d$ are exact (c.f.~\cite{najman2017tightness}). 
This means that the relaxations are exact in Hausdorff metric.
\subsection{Gaussian Probability Density Function}
\label{sec:Relaxations_of_PDF}
The PDF is an auxiliary for the EI acquisition function and is given by $\phi: \mathbb{R} \to \mathbb{R}$ with 
\begin{equation}
\phi(x) \coloneqq \frac{1}{\sqrt{2 \pi}} \cdot \exp\left( \frac{-x^2}{2} \right)	\label{eq:PDF}
\end{equation}
The Gaussian probability density function (PDF) is a nonconvex function for which the McCormick composition rule does not provide its envelopes.
For one-dimensional functions, McCormick~\cite{McCormick.1976} also provides a method to construct envelopes. 
We construct the envelopes of PDF using this method and implement them in our library of intrinsic functions. 
The envelope of the PDF is illustrated in Figure~\ref{fig:McCormickRelaxationsOfPDF} and derived in Appendix~\ref{sec:Derivation_of_relaxations_PDF_Gaussian}.
\begin{figure}[ht]	
	\centering
	\includegraphics[width=0.45\textwidth]{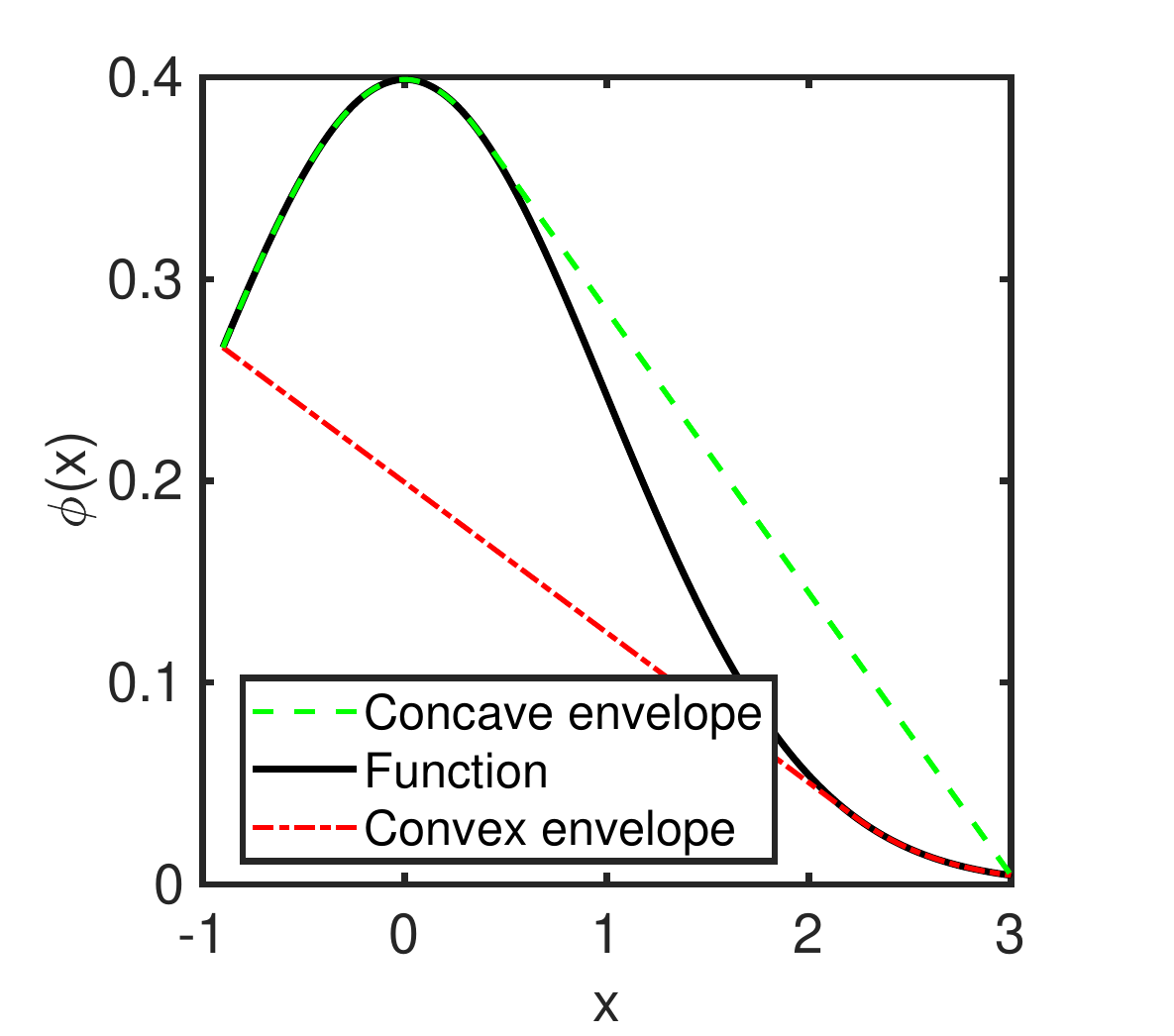}
	\caption{Illustration of the envelope of the Gaussian PDF}
	\label{fig:McCormickRelaxationsOfPDF}  
\end{figure}
\subsection{Gaussian Cumulative Distribution Function}
\label{sec:Relaxations_of_CDF}
The Gaussian cumulative distribution function (CDF) is given by $\Phi: \mathbb{R} \to \mathbb{R}$ with 
\begin{equation}
\Phi(x) \coloneqq \int_{-\infty}^{x} ~\phi(t)~dt = \frac{1+\erf \left( \frac{\sqrt{2} x}{2} \right)}{2}.	\label{eq:CDF}
\end{equation}
The envelopes of the error function are already available in MC\texttt{++} as an intrinsic function and consequently the McCormick technique provides envelopes of the CDF (see Figure~2a in ESI).
In contrast, the error function is not available as an intrinsic function in BARON and a closed-form expression does not exist. 
Thus, a numerical approximation is required for optimization in BARON.
Common numerical approximations of the error function are only valid for $x\ge0$ and use point symmetry of the error function.
To overcome this technical difficulty in BARON, a big-M formulation with additional binary and continuous variables is a possible workaround.
However, this workaround leads to potentially weaker relaxations (see Section~2 in the ESI).
\subsection{Lower Confidence Bound Acquisition Function}
\label{sec:Relaxations_of_LCB}
The lower confidence bound (LCB) (upper confidence bound when considering maximization) is an acquisition function with strong theoretical foundation. 
For instance, a bound on its cumulative regret, i.e., a convergence rate for Bayesian optimization, for relatively mild assumptions on the black-box function is known~\cite{srinivas2009gaussian}.
It is given by $\text{LCB}: \mathbb{R} \times \mathbb{R}_{\ge 0} \to \mathbb{R}$ with 
\begin{equation}
\text{LCB}(\mu, \sigma) \coloneqq \mu - \kappa \cdot \sigma  \notag
\end{equation}
with a parameter $\kappa \in \mathbb{R}_{>0} $.
LCB has not been popular in engineering applications as it requires an additional tuning parameter $\kappa$ and leads to heavy exploration when a rigorous value for $\kappa$ is chosen~\cite{srinivas2009gaussian}.
Recently, LCB has gained more popularity through application as policy in deep reinforcement learning, e.g., by DeepMind~\cite{mnih2013playing}.
LCB is a linear function and thus McCormick relaxations are exact.
\subsection{Probability of Improvement Acquisition Function}
\label{sec:Relaxations_of_PI}
Probability of improvement (PI) computes the probability that a prediction at $x$ is below a given target $f_{\text{min}}$, i.e., $\tilde{\text{PI}}(\boldsymbol{x}) = \mathcal{P}\left( f(\boldsymbol{x}) \le f_{\text{min}} \right)$.
When the underlying function is distributed as a GP with mean $\mu$ and variance $\sigma$, the PI is given by $\text{PI}: \mathbb{R} \times \mathbb{R}_{\ge 0} \to \mathbb{R}$ with 
\begin{equation}
\text{PI}(\mu, \sigma) \coloneqq 
\begin{cases}
\Phi \left( \frac{f_{\text{min}} - \mu}{\sigma} \right) \label{eq:ProbabilityOfImprovement},  &\quad \sigma > 0, \\
0, & \quad  \sigma = 0, ~ f_{\text{min}} \le \mu, \\
1, & \quad  \sigma = 0, ~ f_{\text{min}} > \mu.
\end{cases}
\end{equation}

\begin{figure}[ht]	
	\centering
	\subfloat[]{\centering
		\begin{overpic}[width=0.45\textwidth]{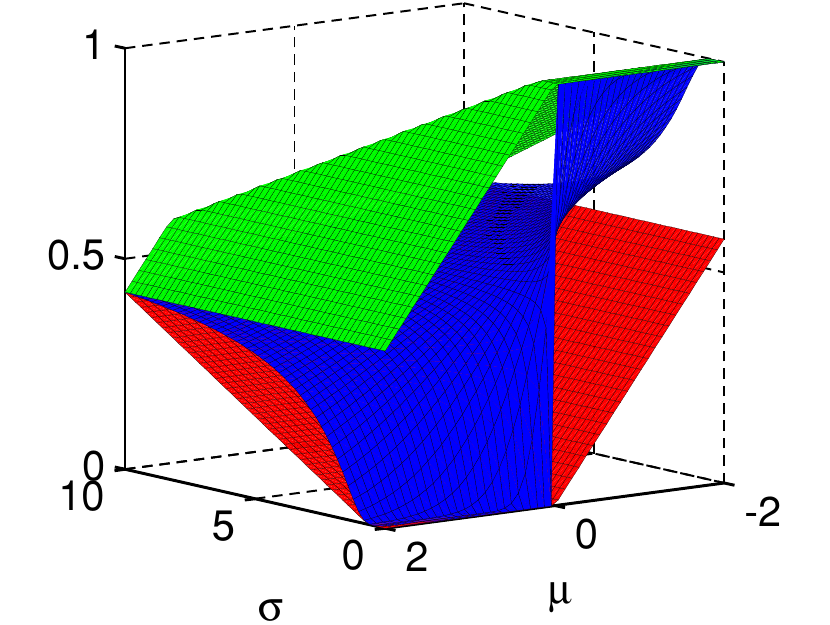}
			\put(76,26){\color{red}Convex}
			\put(76,20){\color{red}relaxation}
			\put(81,55){\color{blue}PI}
			\put(17,66){\color{green}Concave relaxation}
		\end{overpic}\label{fig:PIa}}
	\hfill
	\subfloat[]{\centering
		\includegraphics[width=0.45\textwidth]{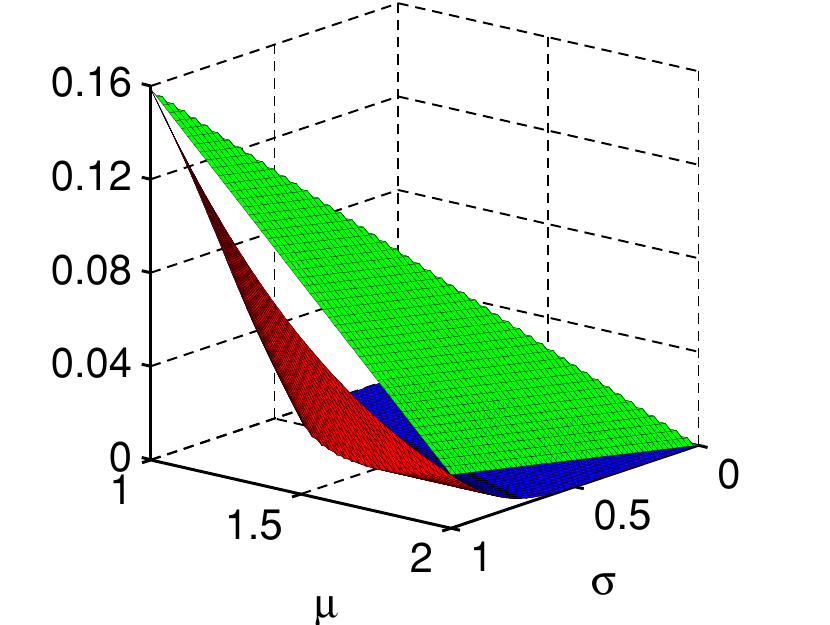}\label{fig:PIb}}
	\caption{Graph of the probability of improvement acquisition function (PI) as in Equation \eqref{eq:ProbabilityOfImprovement} for $f_\textup{min}=0$ along with the developed convex and concave relaxations. (a) On the interval $[-2,2]\times[0,10]$, the relaxations are constructed on the basis of monotonicity properties of PI. (b) On the interval $[1,2]\times[0,1]$, the relaxations are constructed on the basis of componentwise convexity properties via the methods of Meyer and Floudas~\cite{Meyer.2005} and Najman et al. \cite{najman2019convex}}
	\label{fig:PI}  
\end{figure}
Tight relaxations of the PI acquisition function are derived in Section~\ref{sec:Derivation_of_relaxations_PI_acquisition_function} based on a combination of McCormick relaxations, tailored relaxations based on componentwise monononicity, and tailored relaxations based on componentwise convexity/concavity~\cite{Meyer.2005,najman2019convex}. Examples for the resulting relaxations on two different subsets of the domain of $\textup{PI}$ are shown in Figure~\ref{fig:PI}.

\subsection{Expected Improvement Acquisition Function}
\label{sec:Relaxations_of_EI}
EI is the acquisition function that is most commonly used in Bayesian optimization~\cite{jones1998efficient}.
It is defined as $\tilde{\text{EI}}(\boldsymbol{x}) = \E\big[ \max(f_{\text{min}} - f(\boldsymbol{x}) , 0)\big]$.
When the underlying function is distributed as a GP, $\text{EI}: \mathbb{R} \times \mathbb{R}_{\ge 0} \to \mathbb{R}$ is given by 
\begin{equation}
\text{EI}(\mu, \sigma) \coloneqq
\begin{cases}
\left( f_{\text{min}} - \mu \right) \cdot  \Phi \left( \frac{f_{\text{min}} - \mu}{\sigma} \right) + \sigma \cdot \phi \left( \frac{f_{\text{min}} - \mu}{\sigma} \right), & \quad  \sigma > 0 \\
f_{\text{min}} - \mu, & \quad  \sigma = 0, \quad \mu < f_{\text{min}} \\
0 &\quad  \sigma = 0, \quad \mu \ge f_{\text{min}}
\end{cases} 
\label{eq:ExpectedImprovement}
\end{equation}
As noted by Jones et al.~\cite{jones1998efficient}, EI is componentwise monotonic and thus, exact interval bounds can easily be derived.
In Section~\ref{sec:Derivation_of_relaxations_EI_acquisition_function}, we show that EI is convex and we provide its envelopes.
As EI is not available as an intrinsic function in BARON, an algebraic reformulation is necessary that uses Equation~\eqref{eq:ExpectedImprovement} where $\Phi$ is substituted from Equation~\eqref{eq:CDF} with Equation~(1) in ESI and $\phi$ from Equation~\eqref{eq:PDF}.
In addition, some workaround would be necessary for $\sigma = 0$ (e.g., additional binary variable and big-M formulation).
\section{Numerical Results}
\label{sec:Numericalresults}
We now investigate the numerical performance of the proposed method on one core of an Intel Xeon CPU with 2.60~GHz, 128~GB RAM and Windows Server 2016 operating system.
We use MAiNGO version v0.2.1 and BARON v19.12.7 through GAMS v30.2.0 to solve different optimization problems with GPs embedded.
We use SLSQP~\cite{kraft1994algorithm} within the pre-processing of MAiNGO as a local solver. 
The GP models, acquisition functions, and training scripts are available open-source within the MeLOn toolbox~\cite{MeLOn_Git} and the relaxations of the corresponding functions are available through the MC\texttt{++} library used by MAiNGO.
We present three case studies.
First, we illustrate the scaling of the method w.r.t. the number of training data points on a representative test function.
Herein, the estimate of the GP is optimized.
Second, we consider a chemical engineering case study with a chance constraint, which utilizes the variance prediction of a GP.
Third, we optimize an acquisition function that is commonly used in Bayesian optimization on a chemical engineering dataset.
\subsection{Illustrative Example \& Scaling of the Algorithm}
\label{subsec:Numerical Example 1}
In the first illustrative example, the peaks function is learned by GPs.
Then, the GP predictions are optimized on $\tilde{X} = \{x_1, x_2 \in \mathbb{R} \colon -3 \leq x_1, x_2 \leq 3 \}$. 
The peaks function is given by $f_{peaks}: \mathbb{R}^2 \to \mathbb{R}$ with 
\begin{align}
&f_{\text{peaks}}(x_1,x_2) \coloneqq & \notag \\ 
& 3~(1-x_1)^2 \cdot e^{-x_1^2 -~(x_2+1)^2} - 10 \cdot~(\frac{x_1}{5} - x_1^3 - x_2^5) \cdot e^{-x_1^2-x_2^2}- \frac{e^{-(x_1+1)^2 - x_2^2} }{3}     \notag
\end{align}
The two-dimensional function has multiple suboptimal local optima and one unique global minimizer at $\boldsymbol{x}^*\approx[0.228,-1.626]^T$ with $f_{\text{peaks}}(\boldsymbol{x}^*) \approx -6.551$. \newline \indent 
We generate various training data on $\tilde{X}$ using a Latin hypercube sampling of sizes 10, 20, 30, ..., 500.
Then, we train GPs with $k_{\nu=1/2}(d)$, $k_{\nu=3/2}(d)$, $k_{\nu=5/2}(d)$, and $k_{SE}(d)$  covariance functions on the data.
After training, the predictions of the GPs are minimized using the RS and FS formulation to locate an approximation of the minimum of $f_{peaks}$.
We run optimizations in MAiNGO once using the developed envelopes and once using standard McCormick relaxations.
Due to long CPU times, we run optimizations for the FS formulations only for up to 250 data points in MAiNGO.
The whole data generation, training, and optimization procedure is repeated 50 times for each data set.
Thus, we train a total of $10,000$ GPs and run $90,000$ optimization problems in MAiNGO. 
We also solve the FS and RS formulation in BARON by automatically parsing the problem from our C\texttt{++} implementation to GAMS.
This is particularly important in the RS as equations with several thousand characters are generated.
We solve the RS problem for up to 360 and the FS for up to 210 data points in BARON due to the high computational effort. 
The optimality tolerances are set to $\epsilon_{\text{abs. tol.}}= 10^{-3}$ and $\epsilon_{\text{rel. tol.}}= 10^{-3}$ and the maximum CPU time is set to $1,000$ CPU seconds.
The feasibility tolerances are set to $10^{-6}$.
The analysis in this section is based on results for the $k_{\nu=5/2}$ covariance function.
The detailed results for the other covariance functions show qualitatively similar results (c.f. ESI Section~4).  \newline \indent 
In the FS, this problem has $D + 2 \cdot N + 2$ equality constraints and $2 \cdot D + 2 \cdot N + 2$ optimization variables while the RS has $D$ optimization variables and no equality constraints.
Note that for practical applications the number of training data points is usually much larger than the dimension of the inputs, i.e., $N \gg D$. 
The full problem formulation is also provided in ESI Section~3. \newline \indent 
\begin{figure}[h!tb]    
	\centering
	\subfloat[MAiNGO \label{fig:GP_Comparision_of_total_time_Matern_5_MAiNGO}]{\includegraphics[width=0.45\textwidth]{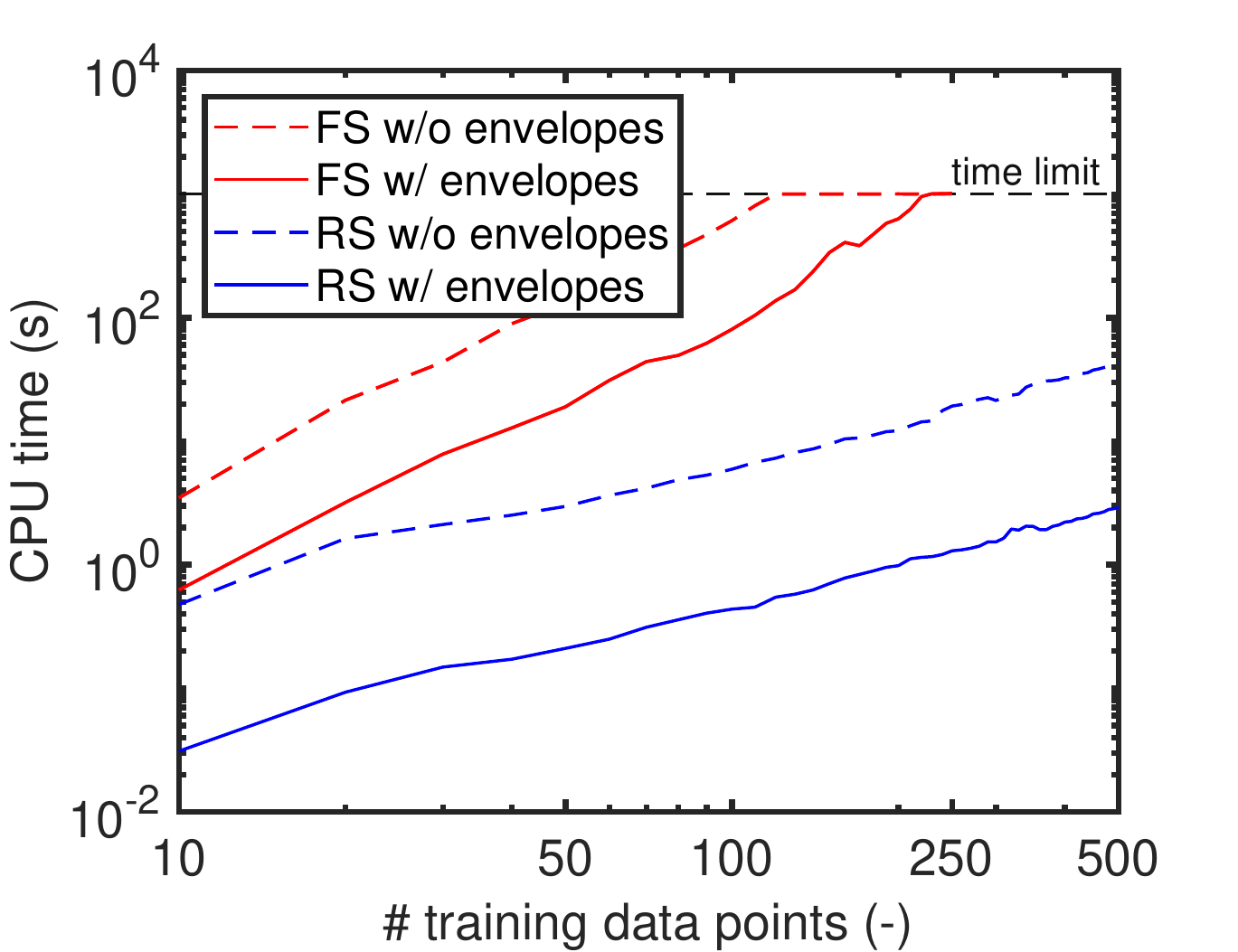}}
	\hfill
	\subfloat[GAMS BARON \label{fig:GP_Comparision_of_total_time_Matern_5_BARON}]{\includegraphics[width=0.45\textwidth]{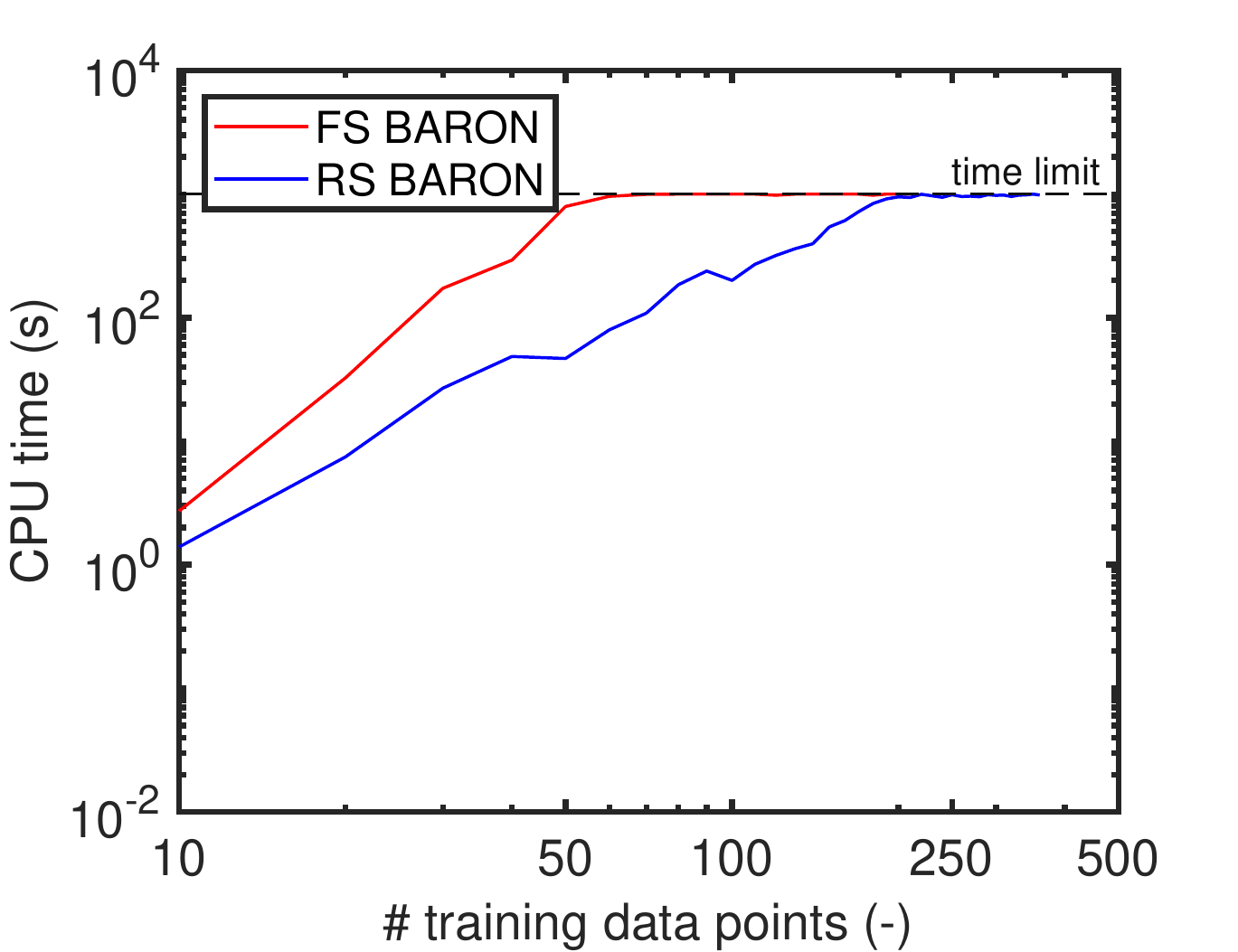}}
	\caption{Comparison of the total CPU time for optimization, i.e., the sum of preprocessing time and B\&B time, of GPs with $k_{\nu=5/2}$ covariance function. The plots show the median of 50 repetitions of data generation, GP training, and optimization. Note that \#points are incremented in steps of 10 and the lines are interpolations between them}
	\label{fig:GP_Comparision_of_total_time_Matern_5}  
\end{figure}
Figure~\ref{fig:GP_Comparision_of_total_time_Matern_5} shows a comparison of the CPU time for optimization of GPs.
For the solver MAiNGO, Figure~\ref{fig:GP_Comparision_of_total_time_Matern_5_MAiNGO} shows that RS formulation outperforms the FS formulation by more than one order of magnitude and shows a more favorable scaling with the number of training data points.
For example, the speedup increases to a factor of 778 for 250 data points. 
Notably, the achieved speedup increases drastically with the number of training data points~(c.f. ESI Section~4).
This is mainly due to the fact that the CPU times for the FS formulations scale approximately cubically with the data points ($\text{CPU}_{FS~w/~env}(N)=1.053\cdot10^{-4} N^{2.958}~\text{sec}$ with $R^2=0.993$) while the ones for the RS scale almost linearly ($\text{CPU}_{RS~w/~env}(N)= 0.0022 \cdot N^{1.156}~\text{sec}$ with $R^2=0.995$). \newline \indent
In general, the number of optimization variables can lead to an exponential growth of the worst-case B\&B iterations and thus runtime. 
In this particular case,
the number of B\&B iterations is very similar for the FS and RS formulation (see Figure~\ref{fig:GP_Branch_and_Bond_Iterations_Analysis_Iterations_MAINGO}).
Instead, for the present problems the number of B\&B iterations is more influenced by the use of tight relaxations.
Figure~\ref{fig:GP_Branch_and_Bond_Iterations_Analysis_CPU_per_iteration_MAiNGO} shows that the CPU time per iteration increases drastically with problem size in the FS while it increases only moderately in the RS.
This indicates that the solution time of the lower bounding, upper bounding, and bound tightening subproblems scales favorably in the RS and that this is the main reason for speedup of the RS formulation in MAiNGO.
This is probably due to the smaller subproblem sizes when using McCormick relaxations in the RS formulation  (c.f. discussion in Section~\ref{sec:Introduction}). \newline \indent
\begin{figure}[h!tb]    
	\centering
	\subfloat[B\&B iterations in MAiNGO \label{fig:GP_Branch_and_Bond_Iterations_Analysis_Iterations_MAINGO}]{\includegraphics[width=0.45\textwidth]{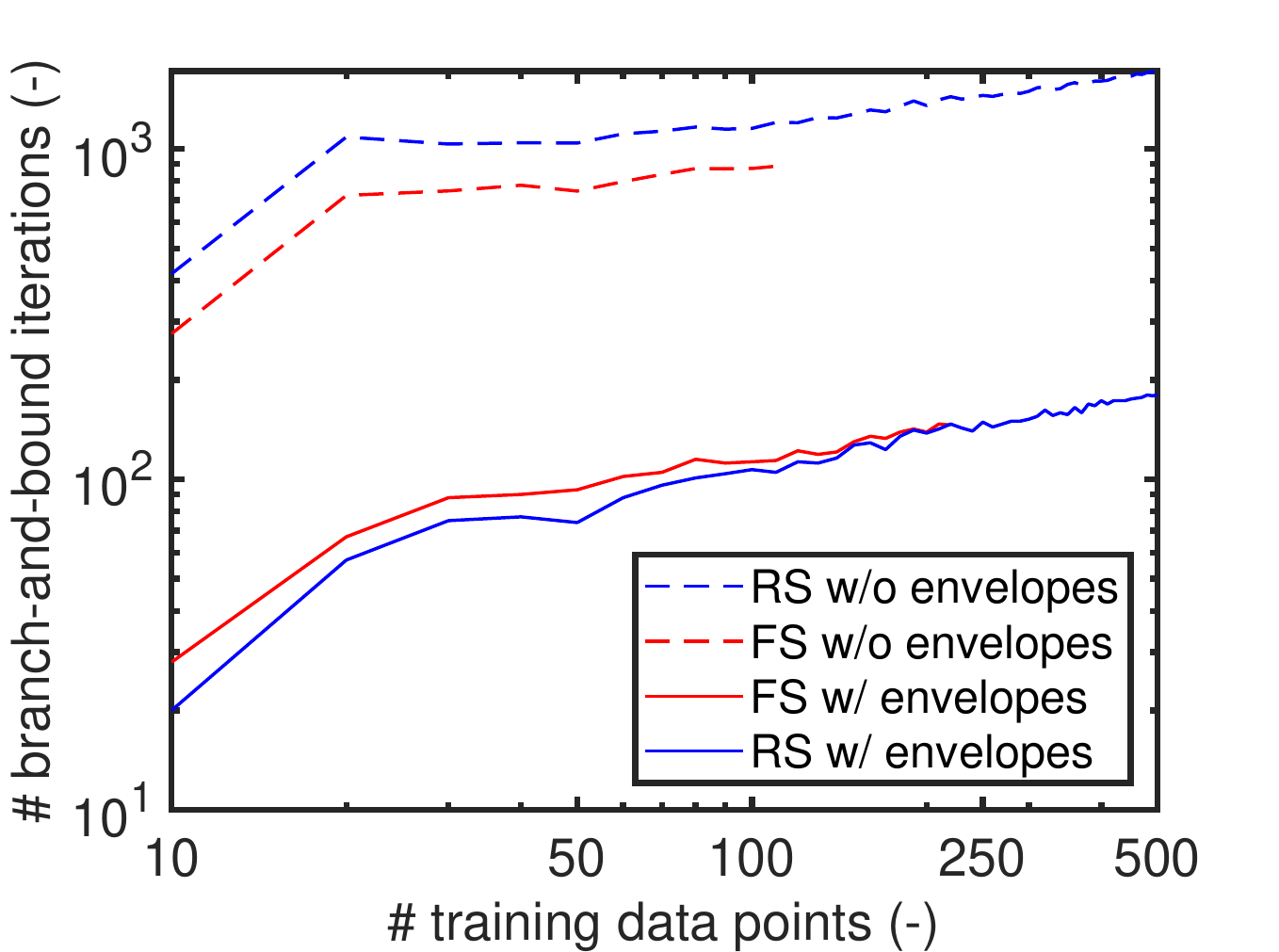}}
	\hfill
	\subfloat[CPU time per B\&B iteration in MAiNGO \label{fig:GP_Branch_and_Bond_Iterations_Analysis_CPU_per_iteration_MAiNGO}]{\includegraphics[width=0.45\textwidth]{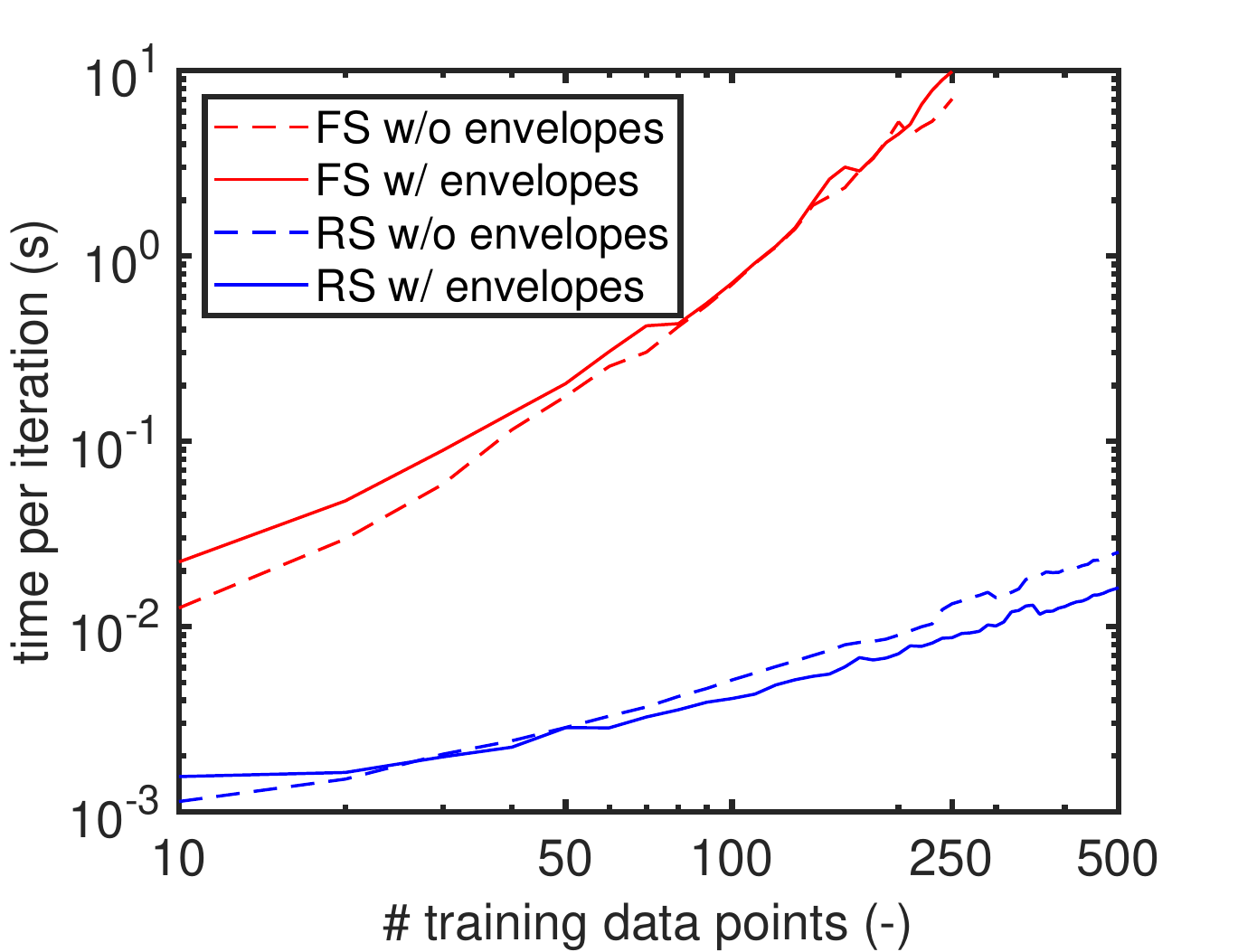}}
	\caption{Comparison of number of B\&B iterations of optimization problems with GPs embedded with $k_{\nu=5/2}$ covariance function. The plots show the median of 50 repetitions of data generation, GP training, and optimization. Note that \#points are incremented in steps of 10 and the lines are interpolations between them}
	\label{fig:GP_Branch_and_Bond_Iterations_Analysis}  
\end{figure}
The use of envelopes of covariance functions also improves computational performance (see Figure~\ref{fig:GP_Comparision_of_total_time_Matern_5_MAiNGO}).
However, this effect is approximately constant over the problem size (c.f. Figure~3 in ESI Section~4).
In other words, the CPU time shows a similar trend for the cases with and without envelopes in Figure~\ref{fig:GP_Comparision_of_total_time_Matern_5_MAiNGO}.
In the RS, the CPU time with envelopes takes on average $7.1\%$ of the CPU time without envelopes ($\approx 14$ times less).
In the FS, the impact of the envelopes is less pronounced, i.e., the CPU time w/ envelopes is on average $15.1\%$ of the CPU time w/o envelopes ($\approx 6.6$ times less).
Figure~\ref{fig:GP_Branch_and_Bond_Iterations_Analysis_Iterations_MAINGO} shows that the envelopes considerably reduce the number of necessary B\&B iterations.
However, the relaxations do not show a significant influence on the CPU time per iteration (see Figure~\ref{fig:GP_Branch_and_Bond_Iterations_Analysis_CPU_per_iteration_MAiNGO}). \newline \indent 
The results of this numerical example show clearly that the development of tight relaxations is more important for the RS formulation than for the FS. 
As shown in Section 3.4.2 of~\cite{bongartz2020_diss}, this effect can be explained by the fact that in RS, it is more likely to have reoccurring nonlinear factors which can cause the McCormick relaxations to become weaker (c.f. also the relationship to the AVM in this case explored in~\cite{Tsoukalas.2014}).
However, in this study, the improvement in relaxations is outweighed by the increase of CPU time per iteration when additional variables are introduced in the FS. \newline \indent
The RS formulation also performs favorably compared to the FS formulation in the solver BARON (see Figure~\ref{fig:GP_Comparision_of_total_time_Matern_5_BARON}).
However, the differences between the CPU times are less pronounced.
In contrast to MAiNGO, the number of B\&B iterations in the FS and RS drastically increase with increasing number of training data points when using BARON (c.f. Figure~4 in ESI).
Also, the time per B\&B iteration is similar between RS and FS.
This is probably due to the AVM method for the construction of relaxations.
The AVM method introduces auxiliary variables for some factorable terms.
Thus, the size of the subproblems in BARON increases with the number of training data points regardless which of the two formulations is used. 
\subsection{Chance-Constrained Programming}
\label{subsec:Numerical Example 2}
Probabilistic constraints are relevant in engineering and science~\cite{charnes1959chance} and GPs have been used in the previous literature to formulate chance constraints, e.g., in model predictive control~\cite{bradford2019stochastic}. \newline \indent
As a second case study, we consider the N-benzylation reaction of $\alpha$ methylbenzylamine with benzylbromide to form desired secondary ($2^\circ$) amine and undesired tertiary ($3^\circ$) amine.
We utilize an experimental data set consisting of 78 data points from a robotic chemical reaction platform~\cite{schweidtmann2018machine}.
We aim to maximize the expected space-time yield of $2^\circ$ amine ($2^\circ$-STY) and ensure that the probability of a product quality constraint satisfaction is above 95\%. 
The $2^\circ$-STY and yield of $3^\circ$ amine impurity ($3^\circ$-Y) are modeled by individual GPs.
Thus, we solve optimization problems with two GPs embedded. 
The chance-constrained optimization problem is formulated as follows
\begin{align}
& \underset{\boldsymbol{x} \in E}{\min} 	&& - \E\big[f_{\text{STY}}(\boldsymbol{x})\big] \notag	\\  
& \text{s.t.} 								&& \mathcal{P}\left( f_{\text{impurity}}(\boldsymbol{x}) \le c \right) \ge 95 \%  \notag
\end{align}
Here, the objective is to minimize the negative of the expected STY. 
This corresponds to minimizing the negative prediction of the GP, i.e., $-m_{\boldsymbol{\mathcal{D}},2^\circ-\text{STY}}$.
The chance constraint ensures that the impurity is below a parameter $c$ with a probability of 95\%.
This corresponds to the constraint $m_{\boldsymbol{\mathcal{D}},3^\circ-Y} +  1.96 \cdot \sqrt k_{\boldsymbol{\mathcal{D}},3^\circ-Y} \le c$ with $c=5$. \newline \indent
The optimization is conducted with respect to four optimization variables: 
(1) the primary ($1^\circ$) amine flow rate of the feed varying between 0.2 and 0.4 mL min$^{-1}$, 
(2) the ratio between benzyl bromide and $1^\circ$ amine varying between 1.0 and 5.0, 
(3) the ratio between solvent and $1^\circ$ amine varying between 0.5 and 1.0, and
(4) the temperature varying between 110 and 150$^\circ$C.  \newline \indent
As this problem is highly multimodal and difficult to solve, we increase the number of local searches in pre-processing in MAiNGO to 500 and increase the maximum CPU time to 24 hours.
The computational performance of the different methods is given in Table~\ref{tab:NumericalResults_GP_ChemicalReaction}.
The results show that none of the considered methods converged to the desired tolerance within the time limit.
The RS formulation in MAiNGO that uses the proposed envelopes outperforms the other formulations and BARON solver as it yields the smallest optimality gap. 
Note that the considered SLSQP solver does not find any valid solution point in the FS in MAiNGO while feasible points are found in the RS.
This demonstrates that the RS formulation can also be advantageous for local solvers.
Note that when using IPOPT~\cite{IPOPT} with 500 multistart points in the FS formulation in MAiNGO, it identifies a local optimum with $f^*=-226.5$ in the pre-processing.
In the ESI, we provide a brief comparison of a few pre-processing settings for this case study. \newline \indent
\begin{table}[ht]
	\centering
	\caption{Numerical results of the N-benzylation reaction optimization with chance constraint~(Subsection~\ref{subsec:Numerical Example 2}). The FS formulation has 330 optimization variables, 326 equality constraints, and 1 inequality constraint. The RS formulation has 4 optimization variables, 0 equality constraints, and 1 inequality constraint. The CPU time limit is 86,400 seconds}
	\label{tab:NumericalResults_GP_ChemicalReaction}
	\begin{tabular}{llrrrrr}
		\hline
		& Solver      & CPU {[}s{]} & Iter.  & UB & LB & Abs. gap \\ \hline
		\multirow{3}{*}{(FS)} 
		& MAiNGO  w/  Env.    	&   86,400   &  26,761  & N/A & -379,2 & N/A  \\
		& MAiNGO  w/o Env.    	&   86,400   &   20,795 & N/A & -1,704.1 &  N/A  \\
		& BARON   				&   86,400   &  219,239 & -226.5   & -53,775.0 &  53,548.5  \\
		\hline
		\multirow{3}{*}{(RS)}  
		& MAiNGO w/  Env.         &  86,400     & $1.6\cdot 10^6$  & -226.5 & -244.8 &  18.3 \\
		& MAiNGO w/o Env.        &  86,400     & $1.1\cdot 10^6$   & -226.5  & -573.2  & 346.7  \\
		& BARON                &  86,400 &  7,927 & -226.5 & -56,394.0 & 56,167.5  \\
		\hline
	\end{tabular}
\end{table}
The best solution of the optimization problem that we found is $x_1 = 0.40$ min$^{-1}$, $x_2 = 1.0$, $x_3 =0.5$, and $x_4 =123.5^\circ$C.
At the optimal point, the predicted $2^\circ$-STY is 226.5 kg m$^{-3}$ h$^-1$ with a variance of 17.1 while the predicted amine impurity is 4.2 \% with a variance of 0.17.
The result shows that the probability constraint ensures a safety margin between the predicted impurity and $c=5$.
Note that the chance constraint is active at the optimal solution point.
\subsection{Bayesian Optimization}
\label{subsec:Numerical Example 3}
In the third case study, we consider the synthesis of layer-by-layer membranes.
Membrane development is a prerequisite for sustainable supply of safe drinking water. 
However, synthesis of membranes is often based on try-and-error leading to extensive experimental efforts, i.e., building and measuring a membrane in the development phase usually takes several weeks per synthesis protocol.
In this case study, we plan to improve the retention of \ce{Na2SO4} salt of a recently developed layer-by-layer nanofiltration membrane system.
The optimization variables are the sodium chloride concentration in the polyelectrolyte solution $c_{NaCl} \in [0,0.5]$ gL$^-1$, the deposited polyelectrolyte mass $m_{PE} \in [0,5]$ gm$^{-2}$, and the number of layers $N_{layer} \in \{1,2,3,..,10\}$. 
The detailed description of the setup is given in the literature~\cite{menne2016precise,rall2019rational}. 
Overall, we utilize 63 existing data points from previous literature~\cite{menne2016precise}.
We identify a promising synthesis protocol based on the EI acquisition function.
Thus, this numerical example corresponds to one step of a Bayesian optimization setup for this experiment.
Global optimization of the acquisition function is particularly relevant due to inherent multimodality of the acquisition functions~\cite{kim2019local} and high cost of experiments.
Note that the experimental validation of this data point is not within the scope of this work. \newline \indent

\begin{table}[ht]
	\centering
	\caption{Numerical results of the membrane synthesis optimization~(Subsection~\ref{subsec:Numerical Example 3}). The FS formulation has 136 optimization variables and 133 equality constraints. The RS formulation has 3 optimization variables}
	\label{tab:NumericalResults_GP_Bayesian_Opt_Membrane}
	\begin{tabular}{llrrrrr}
		\hline
		& Solver      & CPU {[}s{]} & Iter.  & UB & LB & Abs. gap \\ \hline
		\multirow{1}{*}{(FS)} 
		& MAiNGO  w/  Env.   & 3,802 &  25,995   & -2.025 & -2.027 &  0.002      \\
		\multirow{1}{*}{(RS)}  
		& MAiNGO w/  Env.         &  405     & 14,331 & -2.025 & -2.027 & 0.002 \\
		\hline
	\end{tabular}
\end{table}

The computational performance of the proposed method is summarized in Table~\ref{tab:NumericalResults_GP_Bayesian_Opt_Membrane}.
Using the solver MAiNGO, the RS formulation converges approximately 9 times faster to the desired tolerance compared to the FS formulation.
Herein, we use the derived tailored relaxations of the EI acquisition function and envelopes of the covariance functions in both cases.
Notably, the FS requires approximately 1.8 times the number of B\&B iterations compared to the RS formulation, which is much less than the overall speedup.
Thus, the results are in good agreement with the previous examples showing that the reduction of CPU time per iteration in the RS has a major contribution to the overall speedup.
For this example, a comparison to BARON is omitted due to necessary workarounds including several integer variables and function approximations for CDF and EI (c.f., Section~\ref{sec:Relaxations_of_CDF}, \ref{sec:Relaxations_of_EI}). \newline \indent
The optimal solution point of the optimization problem is $c_{NaCl} = 0.362$ gL$^-1$, $m_{PE} =  0$ gm$^{-2}$, and $N_{layer} = 4$. 
The expected retention is 85.32 with a standard deviation $\sigma = 14.8$. 
The expected retention is actually worse than the best retention in the training data of 96.1.
However, Bayesian optimization takes also the high variance of the solution into account, i.e., it is also exploring the space.
EI identifies an optimal trade-off between exploration and exploitation.
\section{Conclusions}
\label{sec:GP_Opt_Conclusionandfuturework}
We propose a RS formulation for the deterministic global solution of problems with trained GPs embedded.
Also, we derive envelopes of common covariance functions or tight relaxations of acquisition functions leading to tight overall problem relaxations. \newline \indent 
The computational performance is demonstrated on illustrative and engineering case studies using our open-source global solver MAiNGO.
The results show that the number of optimization variables and equality constraints are reduced significantly compared to the FS formulation.
In particular, the RS formulation results in smaller subproblems whose size does not scale with the number of training data points  when using McCormick relaxations.
This leads to tractable solution times and overcomes previous computational limitations. 
For example, we archive a speedup factor of 778 for a GP trained on 250 data points.
The GP training methods and models are provided as an open-source module called ``MeLOn - \textbf{M}achin\textbf{e} \textbf{L}earning Models for \textbf{O}ptimizatio\textbf{n}'' toolbox~\cite{MeLOn_Git}. \newline \indent
We thus demonstrate a high potential for future research and industrial applications.
For instance, global optimization of the acquisition function can improve the efficiency of Bayesian optimization in various applications.
It also allows to easily include integer decisions and nonlinear constraints in Bayesian optimization.
Furthermore, the proposed method could be extended to various related concepts such as multi-task GPs~\cite{bonilla2008multi}, deep GPs~\cite{damianou2013deep}, global model-predictive control with dynamic GPs~\cite{wang2006gaussian,bradford2018b}, and Thompson sampling~\cite{chapelle2011empirical,bradford2018efficient}. 
Finally, the proposed work demonstrates that the RS formulation may be advantageous for a wide variety of problems that have a similar structure, including various machine-learning models, model ensembles, Monte-Carlo simulation, and two-stage stochastic programming problems.  \newline \indent

\appendix

\section{Derivations of Convex and Concave Relaxations}
For the sake of simplicity, we use the same symbols in each subsection for the corresponding convex ($F^{cv}$) and concave ($F^{cc}$) relaxations.
To solve the one-dimensional nonlinear equation that arise multiple times in the following, we use Newton's method with 100 iterations and a tolerance of $10^{-9}$.
If this is not successful, we run a golden section search as a backup.
\subsection{Probability Density Function of Gaussian Distribution}
\label{sec:Derivation_of_relaxations_PDF_Gaussian}
In this subsection, the envelopes of the PDF are derived on a compact interval $D = [x^{L}, x^{U}]$. 
As the probability density function is one-dimensional, McCormick~\cite{McCormick.1976} gives a method to construct its envelopes. 
The PDF is convex on $]-\infty,-1]$ and $[1,\infty[$ and it is concave on $[-1,1]$.
Its convex envelope, $F^{cv}: \mathbb{R} \to \mathbb{R}$, and concave envelope, $F^{cc}: \mathbb{R} \to \mathbb{R}$, are given by
\begin{equation}
F^{cv}(x)= 
\begin{cases}
\phi(x),  	  &\quad x^{U} \le - 1, \\
F^{cv}_2(x),  &\quad x^{L} \le -1, \quad -1 \le x^{U} \le 1, \\
\secant(x),   &\quad -1 \le x^{L}, \quad x^{U} \le 1, \\
F^{cv}_4(x),  &\quad -1 \le x^{L}, \quad x^{U} \ge 1, \\
F^{cv}_5(x),  &\quad x^{L} \le -1, \quad x^{U} \ge 1, \\
\phi(x),      &\quad x^{L} \ge 1 \\
\end{cases} \notag
\end{equation}
\begin{equation}
F^{cc}(x)= 
\begin{cases}
\secant(x),  &\quad x^{U} \le - 1, \\
F^{cc}_2(x),  &\quad x^{L} \le -1, \quad -1 \le x^{U} \le 1, \\
\phi(x),      &\quad -1 \le x^{L}, \quad x^{U} \le 1, \\
F^{cc}_4(x),  &\quad -1 \le x^{L} \le 1, \quad x^{U} \ge 1, \\
F^{cc}_5(x),  &\quad x^{L} \le -1, \quad x^{U} \ge 1, \\
\secant(x),   &\quad x^{L} \ge 1 \\
\end{cases} \notag
\end{equation}
where $\secant(x) = \frac{\phi(x^U)-\phi(X^L)}{x^U-x^L} \cdot (x-x^L) + \phi(x^L)$.
$F^{cc}_2: \mathbb{R} \to \mathbb{R}$ is given by:
\begin{equation}
F^{cc}_2(x)= 
\begin{cases}
\frac{\phi(x_{c,2}^{U})-\phi(x^{L})}{x_{c,2}^{U}-x^{L}}\cdot(x-x^{L})+\phi(x^{L}), &\quad x \le x_{c,2}^{U}, \\
\phi(x),  &\quad x > x_{c,2}^{U},
\end{cases} \notag
\end{equation}
where $x_{c,2}^{U} = \min(x^{U*}_{c,2},x^{U})$ and $x^{U*}_{c,2}$ is the solution of $\left.\frac{d\phi}{dx}\right|_x =\frac{\phi(x)-\phi(x^L)}{x-x^L}, \quad x \in [-1,0]$.
$F^{cv}_2: \mathbb{R} \to \mathbb{R}$ is given by:
\begin{equation}
F^{cv}_2(x)= 
\begin{cases}
\phi(x), &\quad x \le x_{c,2}^{L}, \\
\frac{\phi(x^{U}-\phi(x_{c,2}^{L}))}{x^{U}-x_{c,2}^{L}}\cdot(x-x_{c,2}^L)+\phi(x_{c,2}^L),  &\quad x > x_{c,2}^{L},
\end{cases} \notag
\end{equation}
where $x_{c,2}^{L} = \max(x^{L*}_{c,2},x^{L})$ and $x^{L*}_{c,2}$ is the solution of 
$\left.\frac{d\phi}{dx}\right|_x =\frac{\phi(x)-\phi(x^L)}{x-x^L}, \quad x \in [x^L,-1]$.
$F^{cc}_4: \mathbb{R} \to \mathbb{R}$ is given by:
\begin{equation}
F^{cc}_4(x)= 
\begin{cases}
\phi(x),  &\quad x < x_{c,4}^{L},\\
\frac{\phi(x^{U})-\phi(x_{c,4}^{L})}{x^{U}-x_{c,4}^{L}}\cdot(x- x_{c,4}^{L})+\phi(x_{c,4}^{L}), &\quad x \ge x_{c,4}^{L}, 
\end{cases} \notag
\end{equation}
where $x_{c,4}^{L} = \max(x^{L*}_{c,4},x^{L})$ and $x^{L*}_{c,4}$ is the solution of 
$\left.\frac{d\phi}{dx}\right|_x =\frac{\phi(x)-\phi(x^L)}{x-x^L}, \quad x \in [0,x^U]$.
$F^{cv}_4: \mathbb{R} \to \mathbb{R}$ is given by:
\begin{equation}
F^{cv}_4(x)= 
\begin{cases}
\frac{\phi(x_{c,4}^{U})-\phi(x_{L})}{x_{c,4}^{U}-x_{L}}\cdot(x-x_{L})+\phi(x_{L}), &\quad x \le x_{c,4}^{U}, \\
\phi(x),  &\quad x > x_{c,4}^{U},
\end{cases} \notag
\end{equation}
where $x_{c,4}^{U} = \min(x^{U*}_{c,4},x^{U})$ and $x^{U*}_{c,4}$ is the solution of
$\left.\frac{d\phi}{dx}\right|_x =\frac{\phi(x)-\phi(x^L)}{x-x^L}$
on $[x^L,-1]$ if $x^L+x^U<0$ or $[1,x^U]$ if $x^L+x^U>0$
The case $x^L+x^U=0$ is symmetrical and handled separately to avoid numerical issues in Newton.
$F^{cc}_5: \mathbb{R} \to \mathbb{R}$ is given by a combination of $F^{cc}_2$ and $F^{cc}_4$:
\begin{equation}
F^{cc}_5(x)= 
\begin{cases}
\frac{\phi(x_{c,2}^{U})-\phi(x^{L})}{x_{c,2}^{U}-x^{L}}\cdot(x-x^{L})+\phi(x^{L}), &\quad x \le x_{c,2}^{U}, \\
\phi(x),  &\quad x_{c,2}^{U} < x < x_{c,4}^{L}, \\
\frac{\phi(x^{U})-\phi(x_{c,4}^{L})}{x^{U}-x_{c,4}^{L}}\cdot(x- x_{c,4}^{L})+\phi(x_{c,4}^{L}), &\quad x \ge x_{c,4}^{L},
\end{cases} \notag
\end{equation}
$F^{cv}_5: \mathbb{R} \to \mathbb{R}$ is given by:
\begin{equation}
F^{cv}_5(x)= 
\begin{cases}
\frac{\phi(x_{c,5}^{U})-\phi(x^{L})}{x_{c,5}^{U}-x^{L}}\cdot(x-x^{L})+\phi(x^{L}), &\quad x^{L} + x^{U} \ge 0, \quad x \le x_{c,5}^{U}, \\
\phi(x),  &\quad x^{L} + x^{U} \ge 0, \quad x > x_{c,5}^{U}, \\
\phi(x), &\quad x^{L} + x^{U} < 0, \quad x \le x_{c,5}^{L}, \\
\frac{\phi(\phi(x^{U}-x_{c,5}^{L}))}{x^{U}-x_{c,5}^{L}}\cdot(x-x_{c,5}^L)+\phi(x_{c,5}^L),  &\quad x^{L} + x^{U} < 0, \quad x > x_{c,5}^{L},
\end{cases} \notag
\end{equation}
where $x_{c,5}^{U} = \min(x^{U*}_{c,5},x^{U})$ and $x^{U*}_{c,5}$ is the solution of
$\left.\frac{d\phi}{dx}\right|_x =\frac{\phi(x)-\phi(x^L)}{x-x^L}$ on $[x^L,0]$. 
Further, $x_{c,5}^{L} = \max(x^{L*}_{c,5},x^{L})$ and $x^{L*}_{c,5}$ is the solution of
$\left.\frac{d\phi}{dx}\right|_x =\frac{\phi(x)-\phi(x^L)}{x-x^L}$ on $[0,x^U]$. 
\subsection{Probability of Improvement Acquisition Function}
\label{sec:Derivation_of_relaxations_PI_acquisition_function}
In this section, a tight relaxation of the PI acquisition function is derived.
PI is continuous for all $(\mu,\sigma)\in\mathbb{R}\times[0,\infty[~\setminus\{(0,0)\}$, since $\lim_{x\to +\infty}\Phi(x)=1$ and $\lim_{x\to -\infty}\Phi(x)=0$. 
\subsubsection{Monotonicity}
\label{ap:PImonotonicity}
From the gradient of PI on $\mathbb{R}\times(0,\infty)$, 
$$
\label{gradPI}\nabla \text{PI}(\mu,\sigma) = 
-\frac{1}{\sigma^2 \cdot \sqrt{2\pi}} \cdot \exp \left( -\frac{\left(f_\textup{min}-\mu\right)^2}{2 \cdot \sigma^2}    \right) 
\left[
\begin{array}{cc}
\sigma \\
\left(f_\textup{min}-\mu\right)

\end{array}\right],
$$
where $f_\textup{min}$ is a given target, we identify the following monotonicity properties:
\begin{itemize}
\item PI is monotonically decreasing with respect to $\mu$.
\item If $\mu < f_\textup{min}$ then PI is monotonically decreasing with respect to $\sigma$.
\item If $\mu \ge f_\textup{min}$ then PI is monotonically increasing with respect to $\sigma$ (recall that $\textup{PI}(f_\text{min},0)=0$, and $\text{PI}(f_\textup{min},\sigma)=0.5~\forall\sigma\in~]0,\infty[$~).
\end{itemize}
These properties can be used to obtain exact interval bounds on the function values of PI. Furthermore, they can be exploited to construct relaxations as described in Section \ref{ap:PIrelaxations}.

\subsubsection{Componentwise Convexity}
The Hessian of PI on $\mathbb{R}\times(0,\infty[$ is given by
\begin{equation*}
\nabla^2 \text{PI}(\mu,\sigma)=
\left[ \begin {array}{cc} 
-\frac{f_\textup{min}-\mu}{\sigma^3} 			& -\frac{\left(f_\textup{min}-\mu\right)^2-\sigma^2}{\sigma^4} \\
-\frac{\left(f_\textup{min}-\mu\right)^2-\sigma^2}{\sigma^4} 	& \frac{\left(f_\textup{min}-\mu\right) \cdot (\left(2\sigma^2-f_\textup{min}-\mu\right)^2)}{\sigma^5}
\end {array} \right] \cdot \frac{1}{\sqrt{2\pi}} \cdot e^{-\frac{\left(f_\textup{min}-\mu\right)^2}{2\sigma^2}}.
\end{equation*}
The Hessian is indefinite and PI is therefore neither convex nor concave on its whole domain.
However, we find \textit{componentwise convexity} properties on certain parts of the domain, i.e., convexity with respect to one variable when the other is fixed. 
To this end, we divide the domain into the following four sets: 
\begin{itemize}
\item $I_1\coloneqq\{(\mu,\sigma)~|~\mu \leq f_\textup{min}~\land~\mu-f_\textup{min}\geq-\sqrt{2}\sigma\}$,
\item $I_2\coloneqq\{(\mu,\sigma)~|~\mu \geq f_\textup{min}~\land~\mu-f_\textup{min}\leq+\sqrt{2}\sigma\}$,
\item $I_3\coloneqq\{(\mu,\sigma)~|~\mu \leq f_\textup{min}~\land~\mu-f_\textup{min}\leq-\sqrt{2}\sigma\}$,
\item $I_4\coloneqq\{(\mu,\sigma)~|~\mu \geq f_\textup{min}~\land~\mu-f_\textup{min}\geq+\sqrt{2}\sigma\}$.
\end{itemize}
On these sets, PI has the componentwise convexity properties listed in Table~\ref{tab:compConvConcPI}.
\begin{table}[ht]
\centering
\caption{componentwise convexity properties of PI over subsets of its domain}
\label{tab:compConvConcPI}
\begin{tabular}{lll}
	Subset & \multicolumn{2}{c}{componentwise property} \\
	\hline
	$I_1$	&  concave~w.r.t. $\mu$ 			&    convex~~w.r.t. $\sigma$  \\
	$I_2$	&  convex~~w.r.t. $\mu$  			&    concave~w.r.t. $\sigma$  \\
	$I_3$	&  concave~w.r.t. $\mu$  			&    concave~w.r.t. $\sigma$ \\
	$I_4$	&  convex~~w.r.t. $\mu$  			&    convex~~w.r.t. $\sigma$           
\end{tabular}
\end{table}

\subsubsection{Relaxations}
\label{ap:PIrelaxations}
We construct relaxations of PI over a given subset $\mathcal{X}=[\mu^\textup{L},\mu^\textup{U}]\times[\sigma^\textup{L},\sigma^\textup{U}]$ of its domain depending on which of the four sets $I_1$-$I_4$ contains the set $\mathcal{X}$.
If  $\mathcal{X}\subset I_1 \cup I_2$, we use the McCormick relaxations obtained by applying the multivariate composition theorem~\cite{Tsoukalas.2014} to the composition of the rational function $\frac{f_\textup{min}-\mu}{\sigma}$ with $\Phi$ (c.f. Equation~\eqref{eq:ProbabilityOfImprovement}), since these are already very tight. If $\mathcal{X}$ does not fully lie within $I_1 \cup I_2$, the McCormick relaxations get increasingly weaker and we thus resort to other methods as described in the following.

If $\mathcal{X}\subset I_4$, PI is componentwise convex with respect to both variables. 
Therefore, the concave envelope of PI over $\mathcal{X}$ consists of two planes anchored at the four corner points of $\mathcal{X}$ and can be calculated as described by~\cite{Meyer.2005}. 
A tight convex relaxation can be obtained using the method by~\cite{najman2019convex}. 
Since the off-diagonal entries of the Hessian have a constant sign over $I_4$, a sufficient condition for this method is fulfilled (c.f. Corollary 1 in~\cite{najman2019convex}).
An example for the resulting relaxation is shown in Figure~\ref{fig:PIb}. 
Similarly, if $\mathcal{X}\subset I_3$, PI is componentwise concave and we obtain its convex envelope using the method by~\cite{Meyer.2005} and a tight concave relaxation using the method by~\cite{najman2019convex}. \newline \indent
If $\mathcal{X}\subset I_2\cup I_4$, we construct relaxations exploiting the monotonicity properties of PI.
Since for all $(\mu,\sigma)\in I_2 \cup I_4$ we have $\mu \geq f_\textup{min}$, PI is thus monotonically decreasing in $\mu$ and increasing in $\sigma$ over $\mathcal{X}$. 
Therefore, we can construct a convex relaxation $\textup{PI}_{2,4}^\textup{cv}:\mathcal{X}\to[0,1]$ as
\begin{equation}
\textup{PI}_{2,4}^\textup{cv}(\mu,\sigma)\coloneqq\max\left(f^\textup{cv}_{\sigma^\textup{L}}(\mu),f^\textup{cv}_{\mu^\textup{U}}(\sigma)\right),
\label{eq:PIcv24}
\end{equation}
where $f^\textup{cv}_{\sigma^\textup{L}}$ and $f^\textup{cv}_{\mu^\textup{U}}$ are the convex envelopes of the univariate functions
\begin{equation}
f_{\sigma^\textup{L}}:[\mu^\textup{L},\mu^\textup{U}]\to[0,1],\mu\mapsto\textup{PI}(\mu,\sigma^\textup{L}) \notag
\end{equation}
and 
\begin{equation}
f_{\mu^\textup{U}}:[\sigma^\textup{L},\sigma^\textup{U}]\to[0,1],\sigma\mapsto\textup{PI}(\mu^\textup{U},\sigma),
\label{eq:fMuU}
\end{equation}
respectively, i.e., they correspond to the function PI restricted to one-dimensional facets of $\mathcal{X}$ at $\sigma^\textup{L}$ and $\mu^\textup{U}$. Both $f^\textup{cv}_{\sigma^\textup{L}}(\mu)$ and $f^\textup{cv}_{\mu^\textup{U}}(\sigma)$ are valid relaxations of $\textup{PI}$ because of the monotonicity of $\textup{PI}$ over $I_2\cup I_4$. 
By taking the pointwise maximum in~\eqref{eq:PIcv24}, we obtain a tighter relaxation while preserving convexity.
To compute $f^\textup{cv}_{\sigma^\textup{L}}$ and $f^\textup{cv}_{\mu^\textup{U}}$, we can use the method described in Section 4 of \cite{McCormick.1976} because they are one-dimensional functions with a known inflection point. 
To apply this method, we typically need to solve a one-dimensional nonlinear equation, which we do via Newton's method. 
A concave relaxation can be obtained analogously using concave envelopes of PI over one-dimensional facets of $\mathcal{X}$ at $\sigma^\textup{U}$ and $\mu^\textup{L}$.
If $\mathcal{X}\subset I_1\cup I_3$, an analogous method can be used since PI is monotonically increasing in both $\mu$ and $\sigma$. \newline \indent
In the most general case, $\mathcal{X}$ contains parts of all four sets $I_1$-$I_4$. In this case, we can still obtain relaxations by exploiting monotonicity properties. 
In particular, we compute a convex relaxation $\textup{PI}_{1-4}^\textup{cv}:\mathcal{X}\to[0,1]$ as 
\begin{equation}
\textup{PI}_{1-4}^\textup{cv}(\mu,\sigma)\coloneqq\max\left(\tilde{f}^\textup{cv}_{\sigma^\textup{L}}(\mu),f^\textup{cv}_{\mu^\textup{U}}(\sigma)\right),\label{eq:PIcv1-4}
\end{equation}
where $f^\textup{cv}_{\mu^\textup{U}}$ is again the convex relaxation of the univariate function $f_{\mu^\textup{U}}$ as in \eqref{eq:fMuU}, which is still valid because PI is decreasing with respect to $\mu$ on its entire domain. 
In contrast, the convex relaxation of the univariate function at $\sigma^\textup{L}$ as in~\eqref{eq:PIcv24} is not valid because PI is not monotonic with respect to $\sigma$. 
Instead, in~\eqref{eq:PIcv1-4} it is replaced by the convex relaxation $\tilde{f}^\textup{cv}_{\sigma^\textup{L}}$ of the univariate function $\tilde{f}_{\sigma^\textup{L}}:[\mu^\textup{L},\mu^\textup{U}]\to[0,1]$ with 
\begin{equation}
\tilde{f}_{\sigma^\textup{L}}(\mu) \coloneqq
\begin{cases}
\textup{PI}(\mu,\sigma^\textup{L}), &\mu \geq f_\textup{min},\\
\textup{PI}(f_\textup{min},\sigma^\textup{L}) + \frac{\textup{PI}(f_\textup{min},\sigma^\textup{L})-\textup{PI}(\mu^\textup{L},\sigma^\textup{U})}{f_\textup{min}-\mu^\textup{L}}\left(\mu-f_\textup{min}\right), & \text{otherwise}.
\end{cases}
\label{eq:ftilde}
\end{equation}
To see that $\tilde{f}^\textup{cv}_{\sigma^\textup{L}}$ is a valid relaxation of $\textup{PI}$, we first note that by definition it is a relaxation of $\tilde{f}_{\sigma^\textup{L}}$, so it suffices to show that $\tilde{f}_{\sigma^\textup{L}}$ is in turn a relaxation of $\textup{PI}$. The latter is established in the following Lemma.\\
\textbf{Lemma}:
Let \textup{PI} be defined as in \eqref{eq:ProbabilityOfImprovement} and $\tilde{f}_{\sigma^\textup{L}}$ as in \eqref{eq:ftilde}. Then $\tilde{f}_{\sigma^\textup{L}}(\mu)\leq\textup{PI}(\mu,\sigma)~\forall(\mu,\sigma)\in\mathcal{X}\coloneqq[\mu^\textup{L},\mu^\textup{U}]\times[\sigma^\textup{L},\sigma^\textup{U}].$\\
\textbf{Proof}:
Consider first any fixed $\hat{\mu}$ such that $\hat{\mu} \geq f_\textup{min}$. 
In this case, we have $\tilde{f}_{\sigma^\textup{L}}(\hat{\mu})=\textup{PI}(\hat{\mu},\sigma^\textup{L})\leq\textup{PI}(\hat{\mu},\sigma)~\forall\sigma\in[\sigma^\textup{L},\sigma^\textup{U}]$ because of the monotonicity w.r.t $\sigma$ (c.f. Section~\ref{ap:PImonotonicity}). 
Next, consider any $\tilde{\mu}$ such that $\tilde{\mu}<f_\textup{min}$. 
Note that this implies $\mu^\textup{L}<f_\textup{min}$. 
In this case, we have
\begin{align*}
\tilde{f}_{\sigma^\textup{L}}(\tilde{\mu})&=\textup{PI}(f_\textup{min},\sigma^\textup{L}) + \frac{\textup{PI}(f_\textup{min},\sigma^\textup{L})-\textup{PI}(\mu^\textup{L},\sigma^\textup{U})}{f_\textup{min}-\mu^\textup{L}}\left(\tilde{\mu}-f_\textup{min}\right)\\
&= \textup{PI}(f_\textup{min},\sigma^\textup{L})\frac{\tilde{\mu}-\mu^\textup{L}}{f_\textup{min}-\mu^\textup{L}} + \textup{PI}(\mu^\textup{L},\sigma^\textup{U})\frac{f_\textup{min}-\tilde{\mu}}{f_\textup{min}-\mu^\textup{L}}\\
&\leq \textup{PI}(f_\textup{min},\sigma^\textup{U})\frac{\tilde{\mu}-\mu^\textup{L}}{f_\textup{min}-\mu^\textup{L}} + \textup{PI}(\mu^\textup{L},\sigma^\textup{U})\frac{f_\textup{min}-\tilde{\mu}}{f_\textup{min}-\mu^\textup{L}}\\
&\leq \textup{PI}(\tilde{\mu},\sigma^\textup{U})\\
&\leq \textup{PI}(\tilde{\mu},\sigma)~\forall\sigma\in[\sigma^\textup{L},\sigma^\textup{U}],
\end{align*}
where the inequalities follow, in this order, from the monotonicity of PI with respect to $\sigma$ for $\mu\geq f_\textup{min}$, its componentwise concavity with respect to $\mu$ for $\mu<f_\text{min}$, and its monotonicity with respect to $\sigma$ for $\mu< f_\textup{min}$.

\subsection{Expected Improvement Acquisition Function}
\label{sec:Derivation_of_relaxations_EI_acquisition_function}
We now show that the EI acquisition function is convex.
From the Hessian matrix of EI on $\mathbb{R}\times (0,\infty)$
\begin{equation}
\text{Hess}_{EI}(\mu,\sigma) = 
\left[ \begin {array}{cc}
\frac{1}{\sigma} & -\frac{\mu - f_{\text{min}}}{\sigma^2} \\
-\frac{\mu - f_{\text{min}}}{\sigma^2} & \frac{(\mu - f_{\text{min}})^2}{\sigma^3}
\end {array} \right] \cdot \phi\left(-\frac{\mu - f_{\text{min}}}{\sigma}\right), \notag
\end{equation}
we find the eigenvalues $0$ and $\frac{(\mu - f_{\text{min}})^2+\sigma^2}{\sigma^3} \cdot \phi\left(-\frac{\mu - f_{\text{min}}}{\sigma}\right)$.
As $\sigma \ge 0$, EI$(\cdot,\cdot)$ is convex and the envelopes can be constructed directly. \\
\\
\textbf{Contributions of authors}
AMS designed the research concept, ran simulations, and wrote the manuscript.
AMS and DB interpreted the results.
AMS, DB, JN, and AM edited the manuscript.
AMS, XL, DG implemented the GP models and parser.
AMS, DB, JN and TK derived the function relaxations.
JN and DB implemented the function relaxations.
AM is principle investigator and guided the effort.\\
\textbf{Acknowledgements}:
We are grateful to Beno{\^{i}}t Chachuat for providing MC\texttt{++}.
We are grateful to Daniel Menne, Deniz Rall and Matthias Wessling for providing the experimental membrane data for Example 3.
% Finally, we thank the associate editor and the anonymous reviewers for their valuable comments and suggestions.

\section*{Conflict of interest}
The authors declare that they have no conflict of interest.

%\bibliographystyle{spmpsci}
%\bibliography{Bibs}

\end{document}